\documentclass[12pt,reqno]{amsart} 
\usepackage[left=3.5cm,tmargin=2cm,bmargin=2.5cm,centering]{geometry}
\usepackage{siunitx}
\usepackage{hyperref}
\usepackage[utf8]{inputenc}
\usepackage[initials]{amsrefs}  
\usepackage{enumerate}
\usepackage{bm,amsmath,amssymb,amsthm}
\usepackage[all]{xy}
\usepackage{stmaryrd}
\usepackage{amsfonts}
\usepackage{float}
\restylefloat{table}
\theoremstyle{plain}
\newtheorem{theorem}{Theorem}
\newtheorem{conjecture}[theorem]{Conjecture}

\newtheorem{question}[theorem]{Question}
\newtheorem{proposition}[theorem]{Proposition}

\newtheorem*{proposition*}{Proposition}
\newtheorem*{lemma*}{Lemma}
\newtheorem*{corollary*}{Corollary}
\theoremstyle{definition}
\newtheorem*{remark}{Remark}
\newtheorem*{example}{Example}
\newtheorem*{examples*}{Examples}
\newtheorem*{remarks*}{Remarks}
\newtheorem*{remark*}{Remark}

\reversemarginpar

\newcommand{\abs}[1]{\ensuremath{\left|#1\right|}}
\newcommand{\Lquote}[1]{``#1"}
\newcommand{\qtext}[1]{\quad\text{#1}}
\newcommand{\set}[1]{\left\{#1\right\}}
\newcommand{\SB}{\backslash}
\renewcommand{\le}{\leqslant}
\renewcommand{\ge}{\geqslant}
\newcommand{\Mdemi}{\frac{1}{2}}

\newcommand{\cusp}{\mathrm{cusp}}
\newcommand{\A}{\mathbb{A}}
\newcommand{\C}{\mathbb{C}}
\newcommand{\F}{\mathbb{F}}
\newcommand{\N}{\mathbb{N}}
\newcommand{\Q}{\mathbb{Q}}
\newcommand{\R}{\mathbb{R}}
\newcommand{\Z}{\mathbb{Z}}
\newcommand{\BmA}{\mathbb{A}}
\newcommand{\BmC}{\mathbb{C}}
\newcommand{\BmP}{\mathbb{P}}
\newcommand{\BmQ}{\mathbb{Q}}
\newcommand{\BmR}{\mathbb{R}}
\newcommand{\BmZ}{\mathbb{Z}}
\newcommand{\kG}{\mathcal{G}}
\newcommand{\kW}{\mathcal{W}}
\newcommand{\CmO}{\mathcal{O}}
\newcommand{\CmS}{\mathcal{S}}
\newcommand{\FmF}{\mathfrak{F}}
\newcommand{\FmH}{\mathfrak{H}}
\newcommand{\FmS}{\mathfrak{S}}
\DeclareMathOperator {\End}  {\TmE\Tmn\Tmd}
\DeclareMathOperator {\GL} {GL}
\DeclareMathOperator {\Msym} {sym}
\DeclareMathOperator {\Mtr}  {tr}
\DeclareMathOperator {\Sp} {Sp}
\DeclareMathOperator {\SO} {SO}
\DeclareMathOperator {\Msgn}  {sgn}
\def\ra{\rightarrow}
\DeclareMathOperator {\Gal}  {Gal}
\DeclareMathOperator {\PGL} {PGL}
\DeclareMathOperator {\MRe}  {\Re e}
\DeclareMathOperator {\SL} {SL}
\DeclareMathOperator {\SU} {SU}
\DeclareMathOperator {\TLie} {Lie}
\DeclareMathOperator {\GSp} {GSp}
\DeclareMathOperator {\Res}  {res}
\DeclareMathOperator {\Sym} {Sym}

\newcommand{\temp}{\mathrm{temp}}
\newcommand{\Aut}{\mathsf{A}}
\newcommand{\ST}{\mathrm{ST}}
\newcommand{\tO}{\mathrm{O}}
\newcommand{\tU}{\mathrm{U}}
\providecommand{\USp}{\mathrm{USp}}
\providecommand{\Spin}{\mathrm{Spin}}
\newcommand{\even}{\mathrm{even}}
\newcommand{\odd}{\mathrm{odd}}
\newcommand{\disc}{\mathrm{disc}}

\newcommand{\et}{\mathrm{et}}
\newcommand{\Spec}{\mathrm{Spec}\,}
\newcommand{\Rep}{\mathrm{Rep}}
\newcommand{\rec}{\mathrm{rec}}
\newcommand{\FF}{\FmF}
\newcommand{\WD}{\mathrm{WD}}
\renewcommand{\End}{\mathrm{End}}
\providecommand{\kG}{\mathcal{G}}
\providecommand{\kW}{\mathcal{W}}
\providecommand{\geom}{\mathrm{geom}}
\providecommand{\arith}{\mathrm{arith}}
\providecommand{\Frob}{\mathrm{Frob}}

\def\ol{\overline}

\title[]{Families of $L$-Functions and their Symmetry}
\author[P. Sarnak]{Peter Sarnak}
\email{sarnak@math.ias.edu}
\author[S. W. Shin]{Sug Woo Shin}
\email{swshin@math.mit.edu}
\author[N. Templier]{Nicolas Templier}
\email{templier@math.cornell.edu}

\begin{document}

\date{\today}

\begin{abstract}
		In \cite{Sarn:family} the first-named author gave a working definition of a family of automorphic $L$-functions. Since then there have been a number of works \cite{DM06}, \cite{thes:Yang}, \cite{KST:Sp4} \cite{GK:low-lying}, \cite{Kowalski:families} and especially \cite{ST11cf} by the second and third-named authors which make it possible to give a conjectural answer for the symmetry type of a family and in particular the universality class predicted in \cite{KS:bams} for the distribution of the zeros near $s=\frac{1}{2}$. In this note we carry this out after introducing some basic invariants associated to a family.
\end{abstract}
\maketitle
\section{Definition of families and Conjectures}\label{section1}

The zoo of automorphic cusp forms $\pi$ on $G=\GL_n$ over $\Q$ correspond bijectively to their standard completed $L$-functions $\Lambda(s,\pi)$ and they constitute a countable set containing species of different types.
For example there are self-dual forms, ones corresponding to finite Galois representations, to Hasse--Weil zeta functions of varieties defined over $\Q$, to Maass forms, etc.
From a number of points of view (including the nontrivial problem of isolating special forms) one is led to study such $\Lambda(s,\pi)$'s in families in which the $\pi$'s have similar characteristics.
Some applications demand the understanding of the behavior of the $L$-functions as $\pi$ varies over a family. Other applications involve questions about an individual $L$-function.
In practice a \emph{family} is investigated as it arises.

For example the density theorems of Bombieri~\cite{Bombieri:on-largesieve} and Vinogradov~\cite{Vinogradov:density-hypothesis}
are concerned with showing that in a suitable sense most Dirichlet $L$-functions have few violations of the Riemann hypothesis, and as such it is a powerful substitute for the latter.
Other examples are the $\GL_2$ subconvexity results which are proved by deforming the given form in a family (see~\cite{IS00} and~\cite{MV10} for accounts). In the analogous function field setting the notion of a family of zeta functions is well defined, coming from the notion of a family of varieties defined over a base. Here too the power of deforming in a family in order to understand individual members is amply demonstrated in the work of Deligne~\cite{Deligne:WeilII}. In the number field setting there is no formal definition of a family $\FmF$ of $L$-functions.

Our aim is to give a working definition for the formation of a family which will correspond to parametrized subsets of $\Aut(G)$, the set of isobaric automorphic representations on $G(\BmA)$. As far as we can tell these include almost all families that have been studied. For the most part our families can be investigated using the trace formula, monodromy groups in arithmetic geometry and the geometry of numbers, and these lead to a determination of the distribution of the zeros near $s=\frac{1}{2}$ of members of the family. For the high zeros of a given $\Lambda(s,\pi)$, it was shown in \cite{RS96} that the local scaled spacing statistics follows the universal GUE laws (Gaussian Unitary Ensemble). We find that the low-lying zeros (i.e., near $s=\frac{1}{2}$) of a family $\FmF$ follow one of the three universality classes computed in~\cite{book:KS} as the scaling limits of monodromy groups.

For the purpose of defining a family we will assume freely various standard conjectures when convenient. While many of these are well out of reach, important special cases are known and in passing to families they become approachable. We begin by reviewing some notation and invariants associated with individual $\pi$'s.

Any $\pi\in \Aut(G)$ decomposes as an isobaric sum $\pi=\pi_1\boxplus \pi_2\boxplus \cdots \boxplus \pi_r$ with $\pi_j$ an automorphic cusp form on $\GL_{n_j}$, $n_1+n_2+\cdots +n_r=n$ \cite{JPSS83}. Correspondingly $\Lambda(s,\pi)=\Lambda(s,\pi_1)\Lambda(s,\pi_2)\cdots \Lambda(s,\pi_r)$ and this reduces the study to that of cusp forms, which will be our main focus.
Here and elsewhere the central character of $\pi$ is normalized to be unitary and the functional equation relates $\Lambda(s,\pi)$ to $\Lambda(1-s,\tilde \pi)$, where $\tilde \pi$ is the representation contragredient to $\pi$. Furthermore we assume that the central character of $\pi$ is trivial when restricted to $\R_{>0}$ (equivalently is of finite order).
Denote by $\Aut_\cusp(G)$ the subset of cuspidal automorphic representations on $G$. By our normalization this is a countable set.
For $\pi\in \Aut_\cusp(G)$ its conductor $N(\pi)$ is a positive integer defined as the product over appropriate powers of various primes $v$ at which $\pi_v$ is ramified (here $\pi\simeq \otimes_v \pi_v$). It is the integer appearing in the functional equation for $\Lambda(s,\pi)$ (see~\cite{Godement-Jacquet}). The analytic conductor $C(\pi)$ as defined in \cite{IS00} is the product of $N(\pi)$ with a factor coming from $\pi_{\infty}$. The analytic conductor measures the ``complexity'' of $\pi$ (and also the local density of zeros of $\Lambda(s,\pi)$ near $s=\frac{1}{2}$) much like the height of rational points in diophantine analysis. As in that setting the set $S(x)=\set{\pi,\ C(\pi)<x}$ is finite (see~\cite{Brumley:multiplicity-one}). It would be interesting to derive a \emph{Weyl--Schanuel} type theorem for this `universal' family, giving the asymptotic behavior of $S(x)$ as $x$ goes to infinity.\footnote{Brumley--Milicevic~\cite{Brumley-Milicevic:counting} have recently done so for $\GL(2)/\Q$.} We will use $C(\pi)$ to order the elements of a family $\FmF \subset \Aut_{\cusp}(G)$. The root number $\varepsilon(\pi)=\varepsilon(\frac{1}{2},\pi)$ is a complex number of unit modulus that occurs as the sign of the functional equation relating $\Lambda(s,\pi)$ to $\Lambda(1-s,\tilde \pi)$ (\cite{Godement-Jacquet}). We say that $\pi$ is self-dual if $\pi=\tilde \pi$ and in this case $\varepsilon(\pi)=\pm 1$. For a self-dual $\pi$,
$\Lambda(s,\pi\times \pi) = \Lambda(s,\pi,\Msym^2)\Lambda(\pi,s,\wedge^2)$ and $\pi$ is said to be \emph{orthogonal} or \emph{symplectic} according as the first or the second factor above carrying the pole at $s=1$ (in the orthogonal case $\pi$ is a standard functorial transfer of a form on a symplectic group or an even orthogonal group and similarly for the symplectic case from an odd orthogonal group). The symplectic case can only occur if $n$ is even, and if $\pi$ is orthogonal then $\varepsilon(\pi)=1$ (\cite{Lapid:root} and~\cite{Arth:classical}*{Thm 1.5.3.(b)}).

The question of the distribution of $\pi_v$ as $v$ varies over the primes is the generalized \emph{Sato-Tate} problem and its formulation is problematic. Each $\pi_v$ is a point in the unitary dual of $G(\Q_v)$ and according to the generalized Ramanujan conjectures it lies in the tempered dual $\widehat{G(\Q_v)}^{\temp}$ if $\pi$ is cuspidal (see \cite{Sarnak:GRC}). Moreover for $v$ large $\pi_v$ is unramified and hence can be identified with a diagonal unitary matrix $(\alpha_{\pi_v}(1),
\ldots,\alpha_{\pi_v}(n))$ that is a point in an $n$-dimensional torus quotient $T_c/W$, where $T_c$ is the product of $n$ unit circles and $W$ is the permutation group on $n$ letters (we divide by $W$ since the matrix is only determined up to $\GL_n(\C)$ conjugacy). The generalized Sato-Tate conjecture asserts that these $\pi_v$'s become equidistributed with respect to a measure $\mu_{\ST}(\pi)$ on $T_c$ (or more precisely $T_c/W$) as $v\to \infty$.
If $\pi$ corresponds to a finite irreducible Galois representation $\rho$, whose image is denoted $B \subset \GL_n(\C)$, then $\mu_\ST(\pi)$ exists by the Chebotarev density theorem and is equal to the push forward $\mu_B$ of Haar measure on $B$ to the tempered conjugacy classes $G_c^\#\simeq T_c/W$ of $G_c\simeq U(n)$, a maximal compact subgroup which is isomorphic to a compact unitary group.
Langlands \cite{Lang:beyond} suggests that for any $\pi$ there is a (possibly non-connected) reductive algebraic subgroup $B$ of $\GL_n(\C)$ such that $\mu_\ST(\pi)=\mu_{B}$ where the latter denotes the pushforward of the Haar measure on $B\cap G_c$. In~\cite{Serre:NXp} Serre gives a precise formulation in terms of Lie group data and a constructive approach when $\pi$ comes from geometry.
In any case it follows from the analytic properties of $\Lambda(s,\pi)$ and $\Lambda(s,\pi\times \tilde \pi)$ that
\begin{equation}
\begin{aligned}
\int_{T_c} \chi(t) \mu_\ST(t)
&=
\int_{B_c} (\alpha_1(\theta) + \cdots + \alpha_n(\theta)  ) \mu_B(\theta) = 0\\
\int_{T_c} \abs{\chi(t)}^2 \mu_\ST(t)
&=
\int_{B_c} \abs{\alpha_1(\theta) + \cdots + \alpha_n(\theta) }^2 \mu_B(\theta ) =1
\end{aligned}
\end{equation}
where $\chi(t)=\Mtr(t)$.
Hence $B$ is irreducible in $\GL_n(\C)$. In general it may happen that $\mu_{B_1}=\mu_{B_2}$ for $B_1$ not conjugate to $B_2$ in $\GL_n(\C)$ (see \cite{AYY:dim-data}), so that $B$ may not be determined up to conjugacy. For our purposes it is $\mu_B$ that is important, so let $I(T):=I^\#(\GL_n(\C))$ denote the countable set of probability measures that come from irreducible subgroups $B$ as above. Langlands's assertion is that $\mu_\ST(\pi)$ is in $I(T)$ and we will then loosely speak of $\pi$ being of type $B$ if $\mu_\ST(\pi) = \mu_B$, even if $B$ is not unique.

We turn to our formulation of a \emph{parametric family} $\FmF$ of automorphic representations on $G$. $\FmF=(W,F)$ consists of a parameter space $W$ and a map $F:W \to \Aut(G)$, and is based on two very general conjectural means of constructing automorphic forms: spectral and geometric.

\medskip \noindent
\underline{\bf
Harmonic families:}
Let $H$ be a connected reductive algebraic group defined over $\Q$ and $\Aut(H)$ the set of discrete automorphic representations on $H(\A)$. A harmonic (spectral) set $\FmH$ of forms on $H$ is a subset of $\Aut(H)$ consisting of forms $\pi$ which are unramified outside of a finite set of places, or for which $\pi_v\in B_v$ for $v$ in a finite set of places and $B_v$ is a nice subset of positive Plancherel measure in the unitary dual $\widehat{H(\Q_v)}$, or a hybrid of theses conditions.
The important thing is that these sets $\FmH$ can be isolated using the trace formula on $H(\Q) \SB H(\A)$. Let $r:{}^L H\to {}^L G$ be a representation of the corresponding Langlands dual group, then functoriality gives a map $r_*:\FmH \to \Aut(G)$ and defines a parametric family $\FmF=(\FmH,r_*)$ of automorphic representations on $G$.


\medskip \noindent
\underline{\bf
Geometric families:}
These parametric families come from zeta functions which are formed from counting solutions to algebraic equations over finite fields, namely Dedekind zeta functions and Hasse--Weil zeta functions.
Let $\mathcal{W}$ be an open dense subscheme of $\BmA^m_{\Q}=\Spec \Q[W_1,\ldots,W_m]$ (or $\Z[W_1,\ldots,W_m]$ if we work over $\Z$) with $W_1,\ldots,W_m$ transcendental parameters. Let $X$ be a smooth and proper scheme over $\mathcal{W}$ with integral fibers. So specializing the base to $w=(w_1,\ldots,w_m)\in \mathcal{W}(\Q)$ yields a smooth proper variety $X_w$ over $\Q$.

As part of the data defining the corresponding parametric family we restrict the $w$'s locally over $\R$ to lie in a real projective cone $C$ which ensure that the \emph{discriminant} $D(w)$ (see Remark (i) below) corresponding to the family has controlled size in terms of the height of $w$ as a point of $\BmP^m(\Q)$. Put $W=C\cap \mathcal{W}(\Q)$. For the $w$'s in $W$ we get in this way a Hasse-Weil $L$-function (if $X_w$ is zero dimensional, a Dedekind zeta function) on the \'etale cohomology group in a fixed degree $d$ \begin{equation}\label{e:Hasse-Weil}
  L(s, H^d_{\acute{e}t}(X_w\times_\Q \ol{\Q},\ol{\Q}_l))
\end{equation} by specializing to $w$. (See Appendix \ref{sub:Hasse-Weil} below for the definition. It involves a choice of a field isomorphism $\iota:\ol{\Q}_l\simeq \C$, though the expectation is that \eqref{e:Hasse-Weil} is independent of the choice.) Note that the dimension $n$ of the $d$-th cohomology of the closed fibers of $X$ is constant over $W$. Assuming the modularity conjecture (Conjecture \ref{conj:modularity} in Appendix~\ref{sub:Hasse-Weil}) we get a map $F: W\to \Aut(G)=\Aut(\GL_n)$ such that $F(w)$ is the $|\det|^{d/2}$-twist of the automorphic representation corresponding to \eqref{e:Hasse-Weil} (so that $L(s,F(w))=L(s+\frac{d}{2},H^d_{et}(X_w\times_\Q \ol{\Q},\ol{\Q}_l))$). This gives us a parametric family $\FmF=(W,F)$ of automorphic forms.

\begin{remarks*}\begin{enumerate}[(i)]
\item Our aim is a statistical study of members of the family. For parametric families $\FmF=(W,F)$ this means ordering the members according to the sets
\[
\set{w\in W:\ C(F(w))< x},
\]
 and this can be achieved with the caveat that one first replaces $C(F(w))$ by a dominating gauge function $D(w)=\mathrm{Disc}(X_w)$ which  approximates $C(F(w))$. There are many cases for which $F$ is essentially one-to-one and then $F(W)$ is a parametrized subset of $\Aut(G)$. We call such a subset a \emph{parametrized family}, where we can drop the parameter space $W$ since the study of $F(W)$ when ordered by conductor does not depend on the parametrization.

\item Various operations can be performed on parametric families such as union; $\FmF \cup \FmF'$ which is the family with parameter space $W\sqcup W'$ and the corresponding map $F$ or $F'$. For the product $\FmF \times \FmF'$ we take as parameters $W\times W'$ and the map $F(w)\times F(w')$, where the last is the Rankin product giving a form on $GL_{nn'}$ if $\FmF$ is on $\GL_n$ and $\FmF'$ on $\GL_{n'}$. (The product $\pi\times \pi'$ corresponds to the functorial map $\rho\otimes \rho'$ where $\rho$ and $\rho'$ are the standard representations of $\GL_n$ and $\GL_{n'}$.) In this product setting we allow one of the factors to be a singleton in forming the product family. One is tempted to form other boolean operations such as intersections on parametrized  families and this can be done (yielding new families) in many cases. However in general global diophantine equations on the parameters $W$ intervene and these can lead to subsets of $\Aut(G)$ which are not families in our sense and which don't obey any of the predictions below (see Section 3).

\item There are various subsets of $\Aut(G)$ which aren't realized in terms of our general constructions which form natural families and which probably obey the conjectures below. These are defined through Galois and class groups and other arithmetic invariants. For example the set of $\pi$'s which correspond to finite Galois representations, and among these the set of $\pi$'s for which the image of the corresponding Galois representation is a given group $B$ (up to conjugation). Another is the set of Hecke zeta functions of class groups of number fields of a given degree, cf. \S\ref{sub:closing} below. While abelian $H$'s above can be studied to the same extent as our general families using class field theory, we don't know how to study these families in any generality and hence we do not include them as part of the general definition. Note however that one can often produce large parametric subfamilies of these arithmetic `families'.
\item The twist by $|\det|^{d/2}$ in the definition of a geometric family $(W,F)$ is introduced to ensure that the (non-archimedean) local components of $\pi=F(w)$ are unitary, cf. the remark below Conjecture \ref{conj1} and the last paragraph of Appendix \ref{sub:Hasse-Weil}. 
\end{enumerate}
\end{remarks*}

With the definition of a parametric family in place we put forth the basic conjectures about them. These may look far-fetched at first, but unlike the study of individual forms, they can be studied and there is ample evidence (by way of proof) for the conjectures. We will give various examples in Section 2.

For $\pi\in \Aut(G)$ we write the finite part of its standard $L$-function as
\begin{equation}
L(s,\pi) = \prod_{v<\infty} L(s,\pi_v) = \sum_{n=1}^\infty \frac{\lambda_\pi(n)}{n^s}.
\end{equation}
In studying a (harmonic or geometric) parametric family $\FmF=(W,F)$ the first thing one needs to count asymptotically is
\begin{equation} \label{weyllaw}
\abs{\FmF(x)} = \sum_{w:\ C(F(w))<x} 1.
\end{equation}
Since with our normalization there are finitely many automorphic representations $\pi\in \Aut(G)$ of conductor less than $x$, this count is indeed finite as soon as $F$ is finite-to-one.
This means that there are obvious cases that should be excluded, for example if the $L$-map $r$ were to factor through ${}^LH\to W_\Q$ for harmonic families or if $X$ were isotrivial for geometric families. If $F$ is not finite-to-one we impose suitable constraints on the parameter space such as restriction to a projective cone in the geometric setting which renders the finiteness (see \S\ref{sub:n1}).


Also implicit in our definition is the requirement that a family has infinite cardinality.
This infiniteness is not strictly necessary at first since for example Conjecture~\ref{conj1} below reduces to the Sato-Tate conjecture for an individual representation but then as we move on to finer arithmetic invariants and to the universality conjecture this becomes critical.
Thus we assume from now on that the parameter space $W$ is infinite. Then $\abs{\FmF(x)}\to \infty$ as $x\to \infty$ and we expect an asymptotic for $\abs{\FmF(x)}$ that is a power of $x$, possibly with logarithms attached.

The following more general vertical limits should exist as $x\to \infty$ with a modest uniformity in $n\ge 1$:
\begin{equation} \label{quant}
\sum_{w:\ C(F(w))<x} \lambda_{F(w)}(n) = t_{\FmF}(n) \cdot \abs{\FmF(x)}
+ O(n^A \abs{\FmF(x)}^{\delta})
\end{equation}
for some $A<\infty$ and $\delta<1$.
As mentioned before it is understood that in practice the ordering by conductor is often replaced by a closely related ordering involving an approximation in terms of the parameters in the family. Also in some explicit cases one might look at shells $\set{w:\ x<C(F(w))<x+H}$ rather than balls, as smaller sets give finer individual information.

The structure of the limits in~\eqref{quant} can be described in terms of $p$-adic densities. Each $\pi\in \Aut(G)$ determines a point $(\pi_\infty,\pi_2,\pi_3,\pi_5, \ldots)$ in
$\prod_v \widehat{\GL_n(\Q_v)}$, with its product topology (and $\pi_v$ is unramified for $v$ large enough).

\begin{conjecture}[Sato--Tate conjecture for $\FmF$]\label{conj1} There is $p_0=p_0(\FmF)>0$ such that if we order the $w$'s in $W$ by $C(F(w))$ then $F(w)$ is equidistributed in $Y:=\prod\limits_{p\ge p_0} \widehat{\GL_n(\Q_p)}$ with respect to a measure $\mu(\FmF)$ satisfying:

		(i) it is a probability measure and is supported on the tempered spectrum, hence the same holds for $\mu_p(\FmF)$ the projection of $\mu(\FmF)$ on $\widehat{\GL_n(\Q_p)}$,

(ii) it has a decomposition as a convex sum $\mu(\FmF)=\nu_1+\nu_2+\cdots + \nu_r$ of positive measures such that each $\nu_j$ is a product measure on $Y$,

(iii) the average of the $\mu_p(\FmF)$ over $p$ exists and defines the \emph{Sato-Tate} measure $\mu_\ST(\FmF)$ on $T$, that is
\begin{equation}\label{STlimit}
\lim_{x\to \infty}
\frac{1}{x}
\sum_{p_0\le p < x}^{} \log p\cdot \mu_p(\FmF)_{|T} =: \mu_\ST(\FmF)
\end{equation}
(for many families there is no need to average over $p$ as $\lim\limits_{p\to \infty} \mu_p(\FmF)_{|T} = \mu_\ST(\FmF)$),

(iv) $\mu_\ST(\FmF)$ is a probability measure and lies in the convex hull of $I(T)$.
\end{conjecture}
The intuition for (iv) is clear enough, $\mu_\ST(\FmF)$ is a mixture of the measures $\mu_\ST(\pi)$ for the `generic' $\pi$ in $F(W)$. The decomposition asserts that only finitely many $B$-types occur generically in $\FmF$.

\begin{remark}
  A priori $F(w)$ may not define a point in $Y$ but one can simply interpret the equidistribution in the conjecture as asserting in particular that the number of $w$ such that $F(w)$ does not lie in $Y$ is statistically negligible.
  In other words, we need not assume that the local components $\pi_p$ (for $p\ge p_0$) are unitary for each $\pi=F(w)$ to make sense of the conjecture, though we do expect them to be always unitary. For harmonic families, the unitarity of $\pi_p$ is standard (assuming the Langlands functoriality map for $r_*$ is compatible with the transfer of $A$-parameters via $r$) and comes down to the fact that the local $A$-parameters for $\GL_n(\Q_p)$ correspond to unitary representations. For geometric families, the unitarity is known in the case of good reduction but generally conditional on the weight-monodromy conjecture, cf. Remark (iv) above and Appendix \ref{sub:Hasse-Weil}.
\end{remark}

For our purpose only some cruder invariants of $\mu_\ST(\FmF)$ are critical. These are the following indicators:
\begin{equation}\label{i}
\begin{aligned}
i_1(\FmF) &= \int_T \abs{\chi(t)}^2\, \mu_\ST(\FmF)(t) \\
i_2(\FmF) &= \int_T \chi(t)^2\, \mu_\ST(\FmF)(t) \\
i_3(\FmF) &= \int_T \chi(t^2)\, \mu_\ST(\FmF)(t)
\end{aligned}
\end{equation}
where $\chi(t)=\Mtr(t)$.
We note the following equality:
\begin{equation}\label{i3}
i_3(\FmF) = \lim_{x\to \infty} \frac{1}{x}
\sum_{p< x}^{} t_\FmF(p^2) \log p.
\end{equation}

Assuming~\eqref{quant} and the Riemann hypothesis for the relevant $L$-functions one can show the following.
\begin{enumerate}[(i)]
		\item $i_1(\FmF)\ge 1$ and $i_1(\FmF)=1$ iff almost all $F(w)$'s are cuspidal. In this case we say that $\FmF$ is \emph{essentially cuspidal} and for the most part we assume that this is the case. So for our statistical distribution questions the family is in $\Aut_\cusp(G)$.

		\item $0\le i_2(\FmF) \le 1$ and $i_2(\FmF)=1$ iff almost all $F(w)$'s are self-dual and $i_2(\FmF)=0$ iff almost all $F(w)$'s are not self-dual. In the former case we say that $\FmF$ is \emph{essentially self-dual} and in the latter case $\FmF$ is \emph{non self-dual}. Note that $i_2(\FmF)=0 \Rightarrow i_3(\FmF)=0$.

		\item $-1\le i_3(\FmF) \le 1$ and $i_3(\FmF)=1$ iff almost all $F(w)$'s are orthogonal and $i_3(\FmF)=-1$ iff almost all $F(w)$'s are symplectic (called \emph{essentially orthogonal} and \emph{essentially symplectic} respectively).
\end{enumerate}
				The above analysis allows one to compute for any $\FmF$ satisfying~\eqref{quant} the Sato-Tate measures corresponding to the equidistribution of the $F(w)$'s for each of the 3 types. This gives positive measures $\mu_\tU(\FmF)$, $\mu_\tO(\FmF)$ and $\mu_{\Sp}(\FmF)$ on $T$ such that
\begin{equation}
		\mu_\ST(\FmF) = \mu_\tU(\FmF) + \mu_\tO(\FmF) + \mu_{\Sp}(\FmF).
\end{equation}
The proportions of type of $F(w)$ in $\FmF$ are determined from our indicators:
\begin{equation} \label{decomposition}
\begin{aligned}
		\mu_\tU(\FmF)(T)+\mu_\tO(\FmF)(T)+\mu_{\Sp}(\FmF)(T) &=1=i_1(\FmF) \\
		\mu_\tO(\FmF)(T) + \mu_{\Sp}(\FmF)(T) &= i_2(\FmF) \\
\mu_\tO(\FmF)(T) - \mu_{\Sp}(\FmF)(T) &= i_3(\FmF).
\end{aligned}
\end{equation}

As a complement it is helpful to note the following
\[
\int_T \chi(t) \mu_\ST(\FmF)(t) =0,
\]
which follows from the fact that $\FmF$ is essentially cuspidal and hence the absence of pole at $s=1$ for almost all $F(w)$'s. Equivalently the limit
\[
\lim_{x\to \infty}
\frac{1}{x}
\sum_{p< x}
t_{\FmF}(p) \log p
\]
exists and always is equal to zero. This is to be compared with~\eqref{i3} above and~\eqref{rank} below.

The interpretation of these indicators in terms of $B$-types is clear. If $\mu_\ST(\FmF)=\mu_B$ for some $B$ then by classical representation theory of compact groups, $i_1(\FmF)=1$ asserts that $B$ is irreducible in $\GL_n(\C)$, $i_2(\FmF)=1$ asserts that $B$ is self-dual (as a subgroup of $\GL_n(\BmC)$) and $i_3(\FmF)$ is the Frobenius-Schur indicator of $B$ in $\GL_n(\BmC)$.
If assertion (iv) of Conjecture~\ref{conj1} holds, that is $\mu_\ST(\FmF)$ is a convex combination of $\mu_B$'s, then even though this decomposition need not be unique\footnote{Jun~Yu~\cite{Yu:linear-dependence} has given examples of this non-uniqueness.}, collecting the $B$-types according to their indices $i_2$, $i_3$ will reproduce the unique decomposition of $\mu_\ST(\FmF)$ given in~\eqref{decomposition}.

The assertion (ii)  of Conjecture~\ref{conj1} suggests that there is a stronger decomposition $\FmF=\FmF_1 \cup \cdots \cup \FmF_r$, although this is not formally part of the conjecture. Here each subfamily $\FmF_i$ of $\FmF$ has asymptotic density $p_i\in [0,1]$ and $\nu_i=p_i\mu(\FmF_i)$. A family $\FmF_i$ such that $\mu(\FmF_i)$ is a direct product of measures on $\widehat{\GL_n(\Q_p)}$ is irreducible in some sense. For example it is plausible that it implies that its horizontal average $\mu_{\ST}(\FmF_i)$ be of the form $\mu_B$ for some irreducible $B$ as above and thus $\FmF_i$ is essentially homogeneous.

Indeed in many of the examples discussed in Section~\ref{sec:examples} such a $B$ will be shown to exist (see notably~\S\ref{sub:harmonic} and \S\ref{sub:katz}). Then we can say that we have attached a Sato-Tate group $H(\FmF)=B$ to the (irreducible) family $\FmF$.  We abstain from attempting a general conjecture about $H(\FmF)$ for at least two reasons, first because $H(\FmF)$ is not uniquely determined by $\mu_{\ST}(\FmF)$ so that a consistent definition seems hopeless, and second because for certain thin families the existence of $H(\FmF)$ is at the same level of difficulty as the existence of the Langlands group $H_\pi$ for an individual $\pi$ (see \S\ref{sub:twists}).

To put forth our prediction for the distribution of the zeros near $s=\frac{1}{2}$ of members of a family $\FmF$ we need two further invariants attached to the family. The first is the rank, $r(\FmF)$, which is typically zero. The only case where we expect it might not be zero is for geometric families for which $s=\frac{1}{2}$ is a special value of $\Lambda(s,\pi)$ connected with a version of the generalization of the Birch and Swinnerton-Dyer conjecture. In the case of elliptic curves, if there are parametric, global rational solutions to the equations defining $X$ (namely solutions in $\Q(W_1,\ldots,W_m)$) they will specialize to solutions of $X_w$ for $w=(w_1,\ldots,w_m)\in W$. In general one considers not only rational points but rational algebraic cycles as in the conjecture by Tate, Lichtenbaum, Deligne, Bloch--Kato, Beilinson and others.

The rank of the family is concerned with the rate of convergence of $\mu_p(\FmF)$ to $\mu_\ST(\FmF)$, and is defined to be
\begin{equation} \label{rank}
r(\FmF):=  \lim_{x\to \infty} \frac{1}{x} \sum_{p< x} - t_\FmF(p) \sqrt{p} \log p.
\end{equation}
For these geometric families one can show that $t_\FmF(p) \ll p^{-1/2}$, so that~\eqref{rank} measures the next to leading term.

This formula~\eqref{rank} in the context of families and rank of elliptic surfaces has been proposed by Nagao~\cite{Nagao}. For $X$ a family of elliptic curves forming an elliptic surface the equality of $r(\FmF)$ and the rank of $X/\Q(W)$ follows from the Tate conjecture for the surface, see~\cite{RS:rank-elliptic}. The universal distributions for zeros near $s=\frac{1}{2}$ are concerned with fluctuations over the family after removing these persistent zeros at $s=\frac{1}{2}$. In what follows we assume that these have been removed or more simply that $r(\FmF)=0$ (according to definition~\eqref{rank}).\footnote{For  a homogeneous symplectic family of positive rank the third and fourth rows of Conjecture~\ref{conj2} below should read $\epsilon=(-1)^{r(\FmF)}$ and $\epsilon=-(-1)^{r(\FmF)}$ respectively.}

The final invariant of $\FmF$ that we need concerns the symplectic $\pi$'s in $\FmF$. For these the epsilon factor or root number $\varepsilon(\pi)$ can be $+1$ or $-1$ and it is not dictated by the Sato-Tate measure of $\FmF$. According to~\eqref{decomposition} we can decompose the family into essential subfamilies $\FmF_\tU$, $\FmF_\tO$, $\FmF_{\Sp}$ and we would like to decompose $\FmF_{\Sp}$ further as $\FmF_{\Sp,+}$ and $\FmF_{\Sp,-}$ according as $\varepsilon=1$ or $-1$. Since $\varepsilon(\pi)$ is given in terms of a product of local $\varepsilon$-factors at the ramified places of $\pi$, one can compute this decomposition analytically in many cases. However to do so in general involves computing averages over our parametric family of the M\"obius function $\mu$. Namely cancellations in sums
\begin{equation}\label{mobius}
\sum_{w} \mu(M(w))
\end{equation}
where $w$ varies over a large set in $\Z^m$ and $M\in \Z[W_1,\ldots,W_m]$. These are predicted by natural generalization of Chowla's conjectures and are known in special cases~\cite{Helfgott:behaviour-root}.

Assuming these allows one to refine the decomposition~\eqref{decomposition} as
\begin{equation}
\mu_{\ST}(\FmF) =
\mu_{\tU}(\FmF) + \mu_{\tO}(\FmF) + \mu_{\Sp,+}(\FmF) + \mu_{\Sp,-}(\FmF),
\end{equation}
as well as the corresponding decomposition into essentially homogeneous subfamilies. In particular this reduces the study of the distribution of the low-lying zeros (as well as other statistical questions for $\FmF$) to the case of $\FmF$ being one of these four homogeneous families.

We now move to the main statistics of families that we will study, namely low-lying zeros. There are other statistics of interest notably  moments of $L$-values, which are known since the work of Keating--Snaith~\cite{Keating-Snaith} to relate to the symmetry type. See~\cite{CFKRS} and~\cite{cong:park:mich} for a broad review of results and applications (there has been much progress since the appearance of these reviews).
Our definition of families captures most of the examples studied to date (see Section~\ref{sec:examples}), although not all of them (see Section~\ref{sec:3}). Our Conjecture~\ref{conj1} is a precise formulation of all the local statistics expected for families.
In fact our notion of families provides a natural setting for the axiomatic recipies in~\cite{CFKRS}, specifically Conjecture~\ref{conj1} as well as Conjecture~\ref{conj2} below are consistent with the family averaging assumptions made in~\cite{CFKRS}*{p.82}.

Write the zeros of $\Lambda(s,\pi)$ as $\frac12+i\gamma^{(\pi)}_j$ (with multiplicities). For the purpose of studying the zeros near $s=\frac{1}{2}$ we scale the $\gamma^{(\pi)}_j$'s setting
\begin{equation}
\widetilde \gamma^{(\pi)}_j
:=
\gamma^{(\pi)}_j \frac{\log C(\pi)}{2\pi}.
\end{equation}
This normalization is universal (i.e., there are no parameters in this process, the conductor $C(\pi)$ measures the local density).
The four universality classes of distributions determined in~\cite{book:KS} are:
\begin{description}
\item[(1) $\tU(\infty)$] the scaling limit of the distribution near $1$ of eigenvalues of matrices in $U(N)$, $N\to \infty$,
\item[(2) $\Sp(\infty)$] the scaling limit of the distribution near $1$ of eigenvalues of matrices in $\USp(2N)$, $N\to \infty$,
\item[(3) $\SO_{\even}(\infty)$] the scaling limit of the distribution near $1$ of eigenvalues of matrices in $\SO(2N)$, $N\to \infty$,
\item[(4) $\SO_{\odd}(\infty)$] the scaling limit of the distribution near $1$ of the eigenvalues of matrices in $\SO(2N+1)$, $N\to \infty$.
\end{description}
In the theoretical (rather than numerical) study of the $\widetilde \gamma^{(\pi)}_j$'s as $\pi$ varies over $\FmF$ one computes the fluctuation $r$-level densities $W^{(r)}$, $r\ge 1$ (see \cite{book:KS} and also the examples in Section 2), and these determine all other statistics.

We can finally state the
\begin{conjecture}[Universality Conjecture]\label{conj2} Let $\FmF$ be a rank $0$ essentially homogeneous family. Then the low-lying zeros of the members of $\FmF$ follow the laws in the following table:
\begin{table}[H]
		\begin{tabular}[c]{l | l | l}
		Homogeneity Type of $\FmF$ & Symmetry Type of $\FmF$ &  Fluctuation $r$-level density \\
\hline
\quad Non self-dual & \quad $\tU(\infty)$ & \quad $W_0^{(r)}$, $r\ge 1$ \\
\quad Orthogonal & \quad $\Sp(\infty)$  & \quad $W_-^{(r)}$, $r\ge 1$ \\
\quad Symplectic $\varepsilon=1$ & \quad $\SO_{\even}(\infty)$ & \quad $W_+^{(r)}$, $r\ge 1$ \\
\quad Symplectic $\varepsilon=-1$ & \quad $\SO_{\odd}(\infty)$ & \quad $W_-^{(r)}$, $r\ge 1$ \\
\end{tabular}
\end{table}
\end{conjecture}

The $r$-variable densities $W^{(r)}$ are those from~\cite{book:KS}. Note that for the type Symplectic $\varepsilon=-1$, we omit the zero at $s=\frac{1}{2}$, which is there because of the sign of the functional equation when forming the densities of each member. The fact that $W_-^{(r)}$ is entered on lines $2$ and $4$ of this table is surprising but can be related to a similar coincidence  at the level of the Weyl integration formula which is already observed in~\cite{Weyl:classical-groups}.

In the formulation of Conjecture~\ref{conj2} above we have restricted ourselves to homogeneous families. This is for simplicity since one could easily consider families of forms which have mixed types, for example it often happens that essentially symplectic families have a root number that takes both the values $1$ and $-1$ with positive proportion (see Section~\ref{sec:examples} for more examples).
The low-lying zeros of such mixed families will be distributed according to the densities above, with weights determined by the decomposition~\eqref{decomposition}.


The Sato-Tate conjecture for families (Conjecture~\ref{conj1}) is in fact a theorem under mild assumptions as we shall explain with examples in the next section (see notably~\S\ref{sub:katz} for general geometric families and~\S\ref{sub:harmonic} for general harmonic families). The conjecture is independent of the analytic continuation of the corresponding $L$-functions and it only captures a portion of the arithmetic of the families.

This is in contrast to the universality conjecture (Conjecture~\ref{conj2}) which is far reaching. It involves arithmetic cancellations which if true lie much deeper.
Also its formulation relies on the zeros $\gamma_j^{(\pi)}$ and thus assumes the analytic continuation of $\Lambda(s,\pi)$ inside the critical strip, which is often a conditional statement. It seems an interesting question to find a substitute towards an unconditional formulation of the universality conjecture in all cases since the Symmetry Type is an intrinsic invariant of a family that should be independent of functoriality or modularity conjectures.
One important source of additional invariants of families are $p$-adic ones (Selmer groups, $p$-adic $L$-functions, etc) which also can be closely tied with the Symmetry Type, see notably Heath-Brown~\cite{Heath-Brown:Selmer-II}, Bhargava--Shankar~\cite{BS:less-one} as well as the recent~\cite{BKLPR:modeling} and the references there.

Besides theoretical results yielding Conjecture~\ref{conj2} for restricted supports of test functions, an important piece of evidence comes from numerical experiments. There are robust algorithms~\cite{Rubin:computational} to numerically compute the zeros and there is ample and excellent agreement for families of $L$-functions of low degrees.

Another important part of the picture is the function field analogue, where we work with the function field $\F_q(X)$ of a curve $X$ and an $\ell$-adic sheaf $F$ of dimension $d$.
 See \cite{ST11cf}*{p.5} and \cite{Katz:ubiquity} for a discussion.
For example if $\mathcal{F}$ is irreducible self-dual orthogonal then there is a natural pairing on $H^1(X,\mathcal{F})$ which is symplectic invariant by the action of Frobenius. This is consistent with Conjecture~\ref{conj2} and even stronger since it provides a spectral interpretation which is lacking over number fields.

As a corollary to the universality above we conclude that if $n$ is odd, and $\FmF$ a pure self-dual family (i.e., all members are self-dual) then its symmetry type is $\Sp(\infty)$ without any further assumptions (in this case $r(\FmF)=0$ since $s=\frac12$ is not critical in the context of Deligne's special value conjectures~\cite{Deligne:valeurs}; see Appendix \ref{sec:non-criticality}).
Similarly a harmonic family $\FmF$ arising from automorphic forms on split $E_8$, $F_4$ or $G_2$ will have symmetry type $\Sp(\infty)$ since all irreducible representations of their dual groups are self-dual and orthogonal~\cite{book:steinberg}.

\section{Examples}\label{sec:examples}
In this section we collect various examples of families, some old some new, which explicate the notions above and which prove in part the various claims and conjectures. It is this wealth of examples that we have tried to unify.

\subsection{n=1}\label{sub:n1} For $G=\GL(1)$, the set $\Aut(G)$ consists of all the primitive (nontrivial) Dirichlet characters $\chi$ so that parametrized families can be described explicitly. The most basic such family is
\begin{equation}
\FmF^{(2)} = \set{\chi:\ \chi^2=1}.
\end{equation}
In terms of our formation it arises either as all the self-dual forms on $\GL_1$ or as the geometric family coming from the curve $Z^2=W$ over $\Z[W]$, i.e., the Dedekind zeta function of quadratic extensions of $\Q$ after removing the constant factor of $\zeta(s)$. The last gives a parametric family which after a standard square-free sieving argument renders $\FmF^{(2)}$ as a parametrized family. According to Conjecture~\ref{conj2} the Symmetry Type of $\FmF^{(2)}$ should be $\Sp(\infty)$. There is ample evidence for this both numerical and theoretical (see Rubinstein's thesis~\cite{Rubin01}).
In this case where $\GL_1(\C)$ is abelian and $1$-dimensional, $I(T)$ corresponds bijectively to the finite subgroups of $T_c=\set{z:\ \abs{z}=1}$ together with $T_c$ itself. The Sato-Tate measure for $\FmF^{(2)}$ exists and is equal to $\mu_B$ where $B=\set{1,-1}\subset T$. In fact $\mu(\FmF^{(2)})=\prod_v \mu_B$ (that is $\mu_B$ at each place $v$), $r(\FmF)=0$ and $i_1(\FmF^{(2)})=i_2(\FmF^{(2)})=i_3(\FmF^{(2)})=1$.

The precise statement about the low-lying zeros of $L(s,\chi)$ is as follows. For $\chi$ primitive of period $q$ its conductor $N(\chi)$ is $q$ and since $\chi_\infty=1$ or $\Msgn$, the analytic conductor $C(\chi)=q$ as well. To form the $r$-level density sums write the zeros of $\Lambda(s,\chi)$, $\chi\in \FmF^{(2)}$ as
\[
\frac{1}{2} + i \gamma_j^{(\chi)},\qtext{with $j=\pm 1, \pm 2, \cdots$}
\]
where $\gamma_j^{(\chi)}\ge 0$ if $j\ge 1$ and $\gamma_{-j}^{(\chi)}=-\gamma_j^{(\chi)}$.

For $\Phi\in \CmS(\BmR^r)$ even in each variable, form the $r$-level (scaled) densities for the low-lying zeros of $\Lambda(s,\chi)$:
\begin{equation} \label{Dchi}
D(\chi,\Phi):=
\sideset{}{^*}\sum_{j_1,j_2,\cdots,j_r} \Phi\left( \frac{\gamma^{(\chi)}_{j_1}\log C(\chi)}{2\pi},\dots,
\frac{\gamma^{(\chi)}_{j_1}\log C(\chi)}{2\pi}
\right),
\end{equation}
where $*$ denotes the sum is over $j_k=\pm 1,\pm 2,\dots$ and $j_{k_1}\neq j_{k_2}$ if $k_1\neq k_2$. The full $\Sp(\infty)$ conjecture for $\FmF^{(2)}$ is equivalent to
\begin{equation} \label{F2}
\frac{1}{\FmF^{(2)}(x)}
\sum_{\chi\in \FmF^{(2)}(x)}^{}
D(\chi,\Phi)
\rightarrow
\int_{\R^r}^{} \Phi(u) W_-^{(r)}(u)
\mathop{}\!du,\qtext{as $x\to \infty$}
\end{equation}
for any $r\ge 1$ and $\Phi\in \CmS(\R^r)$, where
\[
\begin{aligned}
W_{-}^{(r)}(x_1,\cdots,x_r) &= \det(K_-(x_i,x_j))_{\substack{i=1,\dots,r\\ j=1,\cdots,r}},\\
K_-(x,y) &:= \frac{\sin \pi(x-y)}{\pi(x-y)} - \frac{\sin \pi(x+y)}{\pi(x+y)}
\end{aligned}
\]
and
\begin{equation}
\FmF^{(2)}(x)=
\set{\chi\in \FmF^{(2)}:\
C(\chi) < x
}.
\end{equation}
The first to consider the $1$-level density for this family were \"Ozl\"uk and Snyder~\cite{Ozluk-Snyder}, who proved~\eqref{F2} for $r=1$ and support of the Fourier transform $\widehat \Phi$ of $\Phi$ contained in $(-\frac{2}{3},\frac{2}{3})$. Rubinstein \cite{Rubin01} established~\eqref{F2} for any $r\ge 1$ as long as the support $\widehat \Phi \subset \set{\xi:\ \sum_{j=1}^r \abs{\xi_j}<1}$. Later Gao~\cite{thes:GaoPeng} proved that the limit on the l.h.s. of \eqref{F2} exists for support $\widehat \Phi \subset \set{\xi:\ \sum_{j=1}^r \abs{\xi_j}<2}$
but attempts to prove that his answer agrees with the r.h.s. in \eqref{F2} failed until recently.
What remained was an apparently difficult series of combinatorial identities. These are recently proven in \cite{ERR:low-lying} thus establishing~\eqref{F2} in this bigger range. An interesting feature of their proof is that it uses the function field analogues to verify the identities and in this sense it is similar to the recent proof of the Fundamental Lemma (\cite{Ngo08} and references therein). The point is that replacing $\Q$ by $\F_q(t)$ and computing the analogue $r$-level densities for the family of quadratic extensions of $\F_q(t)$, leads to the same answers and ranges as the case of $\Q$. But now averaging over $q$ and keeping track of uniformity to switch orders leads to the setting in which \cite{book:KS} prove the full $\Sp(\infty)$ conjecture and hence the combinatorial identities must hold in the case of $\Q$! An alternative combinatorial proof of the identities should also be possible along the line of~\cite{CoSnaith:n-correlation}.

\subsection{Number fields and Artin $L$-functions}\label{sub:number-fields}
  The zero dimensional cases of the geometric families are already very rich. Let $K=\Q(W_1,...,W_m)$ with $W_1,...,W_m$ indeterminates and let $f\in K[x]$ be irreducible with splitting field $L$ and Galois group $B$. According to Hilbert's irreducibility theorem the set of $w=(w_1,...,w_m)$ in $\Q^m$ for which $f(x,w)$ is irreducible over $\Q$ and the Galois group of its splitting field $L_w/\Q$ is equal to $B$, is the complement of a thin set (\cite[p.123]{Serre:Mordell-Weil}). We call such $w$'s $f$-generic and these are almost all of the points when counting the $w$'s by height (\cite[\S13.1]{Serre:Mordell-Weil}). Let $\rho:B\ra \GL_n(\C)$ be an irreducible $n$-dimensional representation and let $H=\rho(B)$. To each generic $w$ we have the corresponding irreducible Galois representation $\rho_w:\Gal(L_w/\Q)\ra \GL_n(\C)$. This gives a family of $n$-dimensional Artin $L$-functions $L(s,\rho_w)$ and (conjecturally) automorphic cuspforms $\pi_w$ on $G=\GL_n(\A)$. That is we have a parametrized family $\FmF=(W,F)$ where $F(w)=\pi_w$ for $w$ generic. By the Chebotarev density theorem for each such $w$, the Sato--Tate measure $\mu_{\pi_w}$ exists and is equal to $\mu_H$. So we expect that $\mu_{\ST}(\FmF)=\mu_H$ as well. This is indeed so if we order the $\pi_w$'s by the height of $w$. For $p$ large the asymptotics in \eqref{quant} with $n=p^e$ holds. This follows by considering the $w$'s mod $p$ and then studying the variety $f(x,w_1,...,w_m)=0$ over $\F_p$ and using the theory of Artin's congruence zeta and $L$-functions for curves over finite fields in the case of the variable $w_1$, and~\cites{Weil:riemann,Lang-Weil} in general. This leads to the existence of the vertical limits $\mu_p(\FmF)$ and also that these converge to $\mu_H$ as $p\ra \infty$. That is $\mu_{ST}(\FmF)$ exists and is equal to $\mu_H$ in this ordering. A more appropriate ordering of the $w$'s is by the size of $D(w)$ where $D=D(W_1,...,W_m)$ is the discriminant of $f$. The analogue of \eqref{quant} can be carried out for this ordering as well, at least if $w$ keeps away from directions in which $D(w)$ vanishes. The conductor of $\pi_w$ is essentially the content of $D(w)$ and \eqref{quant} can be carried out if the degree of $D$ is small compared to the number of variables $W$. In all cases we find that $\mu_{\ST}(\FmF)=\mu_H$. Once we have $\mu_H$ the key indicators $i_2(\FmF)$ and $i_3(\FmF)$ (here $i_1(\FmF)=1$) are then determined by the corresponding Schur indicators of $H$. Conjecture \ref{conj2} can be established for $\FmF$ for test functions of limited support (as discussed in \S\ref{sub:n1}) if $D(w)$ is of low degree.

  Some very interesting \emph{parametrized} families arise in connection with Dedekind zeta functions of number fields of fixed degree $k$. For $k=2$ this is the family $\FmF^{(2)}$ in \S\ref{sub:n1}. For $k=3$ consider the parameters $W_1,W_2,W_3,W_4$ and the corresponding binary cubic forms (it is convenient to work projectively here) $f(W)=W_1x^3+W_2 x^2y+W_3 xy^2+W_4y^3$. The Galois group of $f$ over $\Q(W)$ is $S_3$. Let $V(\Q)$ be the $\Q$-vector space of such forms with $w=(w_1,w_2,w_3,w_4)\in \Q^4$. Let $V_{\mathrm{gen}}(\Q)$ denote the points $w\in V(\Q)$ for which the splitting field $L_w$ of $f_w$ is an $S_3$ extension of $\Q$. This together with a fixed irreducible representation $\rho$ of $S_3$ yields a parametric family $\FmF$ as above. The group $\GL_2(\Q)$ acts on $V(\Q)$ by linear change of variables and it preserves the fields $L_w$. The quotient $\GL_2(\Q)\backslash V_{\mathrm{gen}}(\Q)$ parameterizes exactly the $S_3$ splitting fields of degree 3 polynomials over $\Q$ (see~\cite{Wright-Yukie}). In order to count these when ordered by conductor it is best to work over $\Z$ rather than $\Q$ as was done in \cite{DH71} who parametrized and counted the cubic extensions of $\Q$ when ordered by discriminant. With $\GL_2(\Z)$ acting on $V(\Z)$ and $V(\R)$ one determines a fundamental domain $\Omega$ and then orders points in $\Omega(\Z)$ by the discriminant $D(w_1,w_2,w_3,w_4)$ which has degree 4. Furthermore one can sieve to fundamental discriminants and to points in $\Omega_{\mathrm{gen}}(\Z)$. The most delicate point technically is dealing with $w$'s in $\Omega(\Z)$ with $D(w)\le X$ and $w$ near the directions where $D(w)=0$.
  To each $f$ in this parametrized reduced set correspond three conjugate cubic fields $K'_f,K''_f,K'''_f$ gotten by adjoining to $\Q$ one of the roots of $f$ and $\mathrm{disc}(K_f^{(j)})=D(f)$. In this way one obtains a parametrization of the cubic extensions of $\Q$ with Galois group $S_3$. Now $\zeta_{K^{(j)}_f}(s)/\zeta(s)=L(s,\rho_f)$ where $\rho_f$ is the corresponding 2-dimensional irreducible representation of $S_3$. Thus this family $\pi_{\rho_f}$ of $\GL_2$-cuspforms (which are known to exist in this case since $\rho_f$ is dihedral) is the parametrized family $\FmF_3$ of Dedekind zeta functions of cubic extensions. We have $\mu_{\ST}(\FmF)=\mu_H$ where $H$ is the dihedral group $D_3$ in $\GL_2(\C)$. It is orthogonal and hence $i_j(\FmF)=1$ for $j=1,2,3$. In particular $\FmF_3$ has an $\Sp(\infty)$ symmetry. This example is due to Yang \cite{thes:Yang}.

  For $k=4,5$ the parametrization over $\Q$ of degree $k$ extensions with $S_k$ Galois groups in terms of $G(\Q)$ orbits of points in certain $G$-prehomogeneous vector spaces $V$, is given in~\cite{Wright-Yukie}. The theory over $\Z$ as needed to determine the density of such quartic and quintic fields is due to Bhargava (\cite{Bha05}, \cite{Bha10}). In all of these cases (including $k\ge 6$ if they could be suitably parametrized) $\mu_{\ST}(\FmF_k)=\mu_{H_k}$, where $H_k$ is the $k-1$ dimensional representation of $S_k$ realized as the symmetries of the $k-1$ simplex. Since this representation is orthogonal we have $i_j(\FmF_k)=1$, for $j=1,2,3$ and all of these parametrized families have an $\Sp(\infty)$ symmetry. A detailed treatment of families of Artin representations is the subject of~\cite{t:artin}.

\subsection{Families of elliptic curves}\label{sub:elliptic}
We next consider geometric families $E\to \kW$ of curves of genus one.
The $1$-parameter families are geometrically the same as elliptic surfaces fibered over the affine line. The singular fibers are classified by Kodaira and N\'eron and can be determined with Tate's algorithm. A $1$-parameter family is given by polynomials in $\Z[w]$ which are the coefficients of the equation of a plane algebraic curve. A well-studied example is that of quadratic twists of a given elliptic curve which can be written in Weierstrass form as $wy^2=x^3 + ax + b$. It can be viewed as a twist of a fixed elliptic curve with the quadratic family from~\S\ref{sub:n1} (for quadratic twists of any fixed automorphic form see~\S\ref{sub:twists} below). There is a natural 2-parameter family $\FmF^{(\mathrm{ell})}$ given by $y^2=x^3+w_1x +w_2$, where every elliptic curve over $\Q$ appears as a fiber with $a,b\in\Z$.  The discriminant function is $D(w)=4w_1^3+27w_2^2$.  By modularity we obtain in each situation a parametric family $\FmF$ of automorphic cusp forms on $\PGL(2)$.

Conjecture~\ref{conj1} can be verified for each of these families $\FmF$ of elliptic curves and the Sato-Tate measure $\mu_{ST}(\FmF)$ exists with indicators $i_1(\FmF)=i_2(\FmF)=1$ and $i_3(\FmF)=-1$. Hence these families are homogeneous symplectic and correspondingly have symmetry type $O(\infty)$. For $\FmF^{(\mathrm{ell})}$ this follows from a theorem of Birch~\cite{Birch:how} while in general see~\S\ref{sub:katz} below.

There is a caveat that we order the elliptic curves by height rather than conductor. Ordering by height for $\FmF^{(\mathrm{ell})}$ means that we restrict to a box, 
\[
\max(4|w_1|^3,27|w_2|^2)<x, \quad \text{with $x\to \infty$.}
\] 
It is desirable to be able to order by conductor $C(w)<x$ with $x\to \infty$ which yields interesting questions related to the square-free sieve for the discriminant polynomial $D(w)$. For $\FmF^{(\mathrm{ell})}$ it follows from~\cites{FNT:Szpiro,Duke-Kowalski} that the number of non-isogeneous semistable elliptic curves of conductor $C(w)<x$ is at least $x^{\frac56}$ and at most $x^{1+\varepsilon}$. The average conductor is also important and it leads one to consider the ratio $\frac{\log C(w)}{\log |D(w)|}$ which is less than $1$ and according to a conjecture of Szpiro should be greater than $\frac16-\varepsilon$ with a finite number of exceptions. For $\FmF^{(\mathrm{ell})}$ the ratio can be shown to be one on average using the square-free sieve which is known for polynomials in $2$-variables of degree $\le 6$ by Greaves~\cite{Greaves:power-free} (for $1$-parameter families it is known for degree $\le 3$ by~\cite{Hooley:book-sieve}).

The next interesting invariant is the rank $r(\FmF)$ defined in~\eqref{rank}. For $\FmF^{(\mathrm{ell})}$ it follows from~\cite{Birch:how} that $t_{\FmF^{(\mathrm{ell})}}(p)\ll p^{-1}$ and thus $r(\FmF^{(\mathrm{ell})})=0$. For a $1$-parameter family it is shown in Miller~\cite{Miller:rational-families} using the Tate conjecture proven in Rosen-Silverman~\cite{RS:rank-elliptic} that $r(\FmF)$ coincides with the rank of the elliptic surface over $\Q(w)$. There are examples of $1$-parameter families where $r(\FmF)$ is greater than $18$ and indeed such families have been used via specialization to produce rational elliptic curves of high rank~\cite{Elkies:high-rank}.

Mazur showed that there are finitely many possibilities for the torsion subgroup of elliptic curves over $\Q$. Harron-Snowden~\cite{Harron-Snowden} recently established various bounds towards counting elliptic curves with prescribed torsion subgroup. In the process they actually show that for each prescribed torsion subgroup, elliptic curves are parametrized by a corresponding moduli space which is close to being an open subscheme of the affine line $\A^1$. Thus these are parametric families according to our definition (e.g., see~\cite{Harron-Snowden}*{\S3} where each family is explicitely given by polynomial equations with one free parameter). 

The root number is the subtlest of the invariants. In the family $(7+7w^4)y^2=x^3-x$ found by Cassels--Schinzel~\cite{CS:families}, the root number $\epsilon(\frac12,E_w)=-1$ for all $w\in \Z$,\footnote{In fact $\epsilon(\frac12,E_w)=-1$ also if we let $w\in \Q$ which should be viewed a $2$-parameter family by writing $w=\frac{w_1}{w_2}$ and ordering by height $\max(|w_1|,|w_2|)<x$.} whereas the rank $r(\FmF)=0$. Another example~\cite{Washington:simplest-cubic} is the $1$-parameter family $y^2=x^3+wx^2-(w+3)x+1$ which has root number $\epsilon(\frac12,E_w)=-1$ for all $w\in \Z$ and for which $r(\FmF)=1$. Thus the rank $r(\FmF)$ and the root numbers of member of $\FmF$ can behave independently from one another and this explains why in Section~\ref{section1} we treat them as distinct invariants.

The average root number is governed by the polynomial $M\in \Z[w_1,\ldots,w_m]$ whose zero set is the locus of the fibers $E_w$ with nodal (multiplicative) singularity. Note that $M$ is a polynomial factor of the discriminant $D$. It is shown by Helfgott~\cite{Helfgott:behaviour-root} how the average root number in these cases is reduced to sums of the type~\eqref{mobius} and thus if $M$ is non-constant, that is if the family has at least one nodal geometric fiber, then the average root number should be zero. In the two examples from the preceding paragraph $M$ is constant and indeed one can find in~\cites{Rizzo:average-root,Helfgott:behaviour-root} further examples of families of elliptic curves with $M$ constant, where the average root number can assume any value in a dense subset of $[-1,1]$. The sum~\eqref{mobius} can be estimated unconditionally for polynomials of low degree, for example~\cite{Helfgott:parity-irreducible}
\begin{equation}\label{helfgott}
 \frac{1}{x^2} \sum_{|w_1|,|w_2|<x} \mu(w_1^3+2w_2^3) = o(1),\quad
\text{as $x\to \infty$.}
\end{equation}
 An example where the root number is shown  to average to zero unconditionally is the $2$-parameter family $y^2=x(x+w_1)(x+w_2)$ which contains every elliptic curve over $\Q$ with full rational $2$-torsion $(\Z/2\Z)^2$ as a fiber over $(w_1,w_2)\in \Z^2$.  The case of $\FmF^{(\mathrm{ell})}$ is more difficult. The method of proof of~\eqref{helfgott} is closely related to the work of Friedlander--Iwaniec and Heath-Brown on primes represented by polynomials in $2$-variables.

The upshot is that Conjecture~\ref{conj2} is verified for families of quadratic twists in~\cite{Rubin01}, for $\FmF^{(\mathrm{ell})}$ in~\cites{You06,Baier-Zhao:low-lying} and under the above assumptions for $1$-parameter families in~\cite{Miller:rational-families}. This yields upper-bounds for the average analytic\footnote{The average rank of Selmer groups, which yields upper-bounds for the average Mordell-Weil rank, can be bounded by other methods, see~\cite{Heath-Brown:Selmer-II,FIMR:spin} for the $1$-parameter families of quadratic twists and~\cite{BS:less-one,Bhargava-Skinner} for $\FmF^{(\mathrm{ell})}$.} rank as a corollary, see for example the articles in the proceedings~\cite{book:ranks-rdm}.
\subsection{Dwork families}\label{sub:dwork}
  In this section we investigate a certain parametric family of Dwork hypersurfaces, which were prominent examples in Dwork's detailed study of hypersurfaces in 1960's. (See introduction of \cite{Kat09} for a commentary on the literature.) Let $U=\Spec \Z[\frac{1}{n+1},w]$, a subscheme of the affine line over $\Z[\frac{1}{n+1}]$. Consider the subscheme $X$ of $\BmP^{n}_{U}$ cut out by the equation
  $$\sum_{i=0}^n x_i^{n+1} = (n+1)w \prod_{i=0}^n x_i,$$
  where $(x_0:\cdots: x_n)$ and $w$ are the coordinates for $\BmP^n$ and $U$, respectively. The family $X\ra U$ is a family of elliptic curves for $n=2$ and that of K3 surfaces for $n=3$. In general the fibers of $X\ra U$ have dimension $n-1$, so the cohomology in degree $n-1$ is the most interesting. We could work with the whole of $H^{n-1}$ cohomology but it is convenient to deal with a piece of cohomology by exploiting a group action on $X$. Let $\mu_{n+1}$ be the set of $(n+1)$-st roots of unity. (One may view $\mu_{n+1}$ as a group scheme over $U$.) Let $H$ be the quotient group $(\mu_{n+1})^{n+1}/\Delta(\mu_{n+1})$, where $\Delta$ is the diagonal embedding. Then $H$ acts on $X$ by letting $(\alpha_0:\cdots:\alpha_n)$ act by $(x_0:\cdots:x_n)\mapsto (\alpha_0x_0:\cdots : \alpha_n x_n)$ on $X$. Let $H_0$ denote the subgroup of $H$ which is a quotient of $\{(\alpha_0:\cdots:\alpha_n): \prod_{i=0}^n \alpha_i=1\}$ by $\Delta(\mu_{n+1})$.

         Consider the setup and notation for geometric families in Section~\ref{section1}. Take $C$ to be the set of $w\in \Z$ such that $w\nmid (n+1)$, viewed as a set of closed points of $U$. Denote by $X_w$ the fiber of $X$ over $w\in U$. Use the discriminant function $D(w)=w^{n+1}-1$ on $C$. Define the map $F:C\ra \Aut(GL_n)$ such that $F(w)$ is the $|\det|^{\frac{n-1}{2}}$-twist of the automorphic representation corresponding to the $\Gal(\ol{\Q}/\Q)$-representation
    \begin{equation}\label{e:cohomology-Dwork}
      H^{n-1}_{\et}(X_w\times_\Q \ol{\Q},\Q_l)^{H_0}
    \end{equation}
    via Conjecture \ref{conj:modularity} (or Conjecture \ref{conj:correspondence}). Note that $F(w)\in \Aut(GL_n)$ since \eqref{e:cohomology-Dwork} has dimension $n$ over $\Q_l$ as can be shown by computing its dimension for $w=0$ (\cite[Lem 1.1]{DMOS82}, cf. \cite[Lem 1.1]{HSBT10}). Since $X_w$ has good reduction modulo $p$ whenever $p\nmid D(w)$, cf. \cite[\S3]{Kat09}, the Galois representation \eqref{e:cohomology-Dwork} is unramified at such $p$, hence $F(w)$ should be unramified outside the prime divisors of $D(w)$.

Suppose that $n$ is \emph{even}. 
The monodromy of the Dwork family $\FmF$ is shown by Dwork to be the full symplectic group (if one is only interested in the symplectic pairing it can also be constructed by Poincar\'e duality, cf. \cite[Lem 1.10, Cor 1.11]{HSBT10}).
    The two main conjectures from Section~\ref{section1} yield the following: first, $\mu_{\ST}(\mathfrak{F})$ arises from the push-forward of a Haar measure on a maximal compact subgroup of $\Sp(n,\C)$ in $\GL(n,\C)$. This is proved as explained in~\S\ref{sub:katz} below using the Deligne--Katz equidistribution theorem. In other words the family has a Sato-Tate group $H(\FmF) = \Sp(n,\C)$.
Second, Conjecture~\ref{conj2} says that the Symmetry Type of $\mathfrak{F}$ should be a superposition of $\SO_{\even}(\infty)$ and $\SO_{\odd}(\infty)$. The superposition depends on the distribution of $\varepsilon=1$ and $\varepsilon=-1$ which we expect will be $\num{50}\%$.

    Finally when $n$ is \emph{odd}, \eqref{e:cohomology-Dwork} is even dimensional and equipped with a perfect symmetric pairing and the exact monodromy is also computed by Dwork. Thus in this case $\mu_{\ST}(\mathfrak{F})$ arises from an even orthogonal group and Conjecture~\ref{conj2} says that the Symmetry Type of $\mathfrak{F}$ should be $\Sp(\infty)$. It would be desirable to test all these low-lying zeros predictions for this family numerically.



\subsection{Harmonic families and Plancherel equidistribution}\label{sub:harmonic}
Consider a spectral set $\FmH\subset \Aut(H)$ of automorphic representations of a connected reductive group $H$ over $\Q$ and an $L$-map $r:{}^LH\to {}^L\GL_n$. These data give rise to a harmonic family $\mathfrak{F}$. We discuss the Sato-Tate equidistribution for $\mathfrak{F}$ as formulated in Conjecture~\ref{conj1}. In fact we need not assume the functoriality conjecture for $r$ to make sense of the conjecture. Namely for each $\sigma\in \FmH$ unramified outside of the finite set of places $S$, we can attach~\cite{Borel:l-fn} the partial $L$-function $L^S(s,\sigma,r)$, which should be the partial $L$-function for $r_*\sigma$ if we assumed that $r_*\pi$ was an automorphic representation of $GL_n$. The prime $p_0$ is chosen large enough so that $p\ge p_0 \Rightarrow p\not\in S$ and thus the unramified representation $\pi_p = r_*\sigma_p$ is known.

The asymptotic count of~\eqref{weyllaw} is a Weyl's law or limit multiplicity problem. This has a long history with a vast literature. For limit multiplicities for towers of subgroups it starts with the classical article of DeGeorge--Wallach~\cite{DeGeor-Wall}.
In the case that $\sigma\in \FmH$ have discrete series $\sigma_\infty$ at infinity the asymptotic count is well-understood and it is natural to first focus on this case for studying harmonic families. See the end of this subsection for a discussion of the Maass forms case.

The next step is the quantitative equidistribution~\eqref{quant} which is much more difficult to obtain.
The $\PGL(2)$ case is treated in~\cite{ILS00}, see \S\ref{sub:gl2} below.
A generalization to higher rank groups was recently achieved by the second and third-named authors~\cite{ST11cf}.  

To fix notation, the spectral set $\FmH$ will be the set of automorphic representations $\sigma$ of $H(\BmA)$ which are cohomological at infinity with regular weight. (This means that the infinite component of $\sigma$ has nonzero Lie algebra cohomology against an irreducible algebraic representation with regular highest weight.) Such $\sigma$ is always cuspidal by a theorem of Wallach. If we consider the weight aspect it will be convenient to fix a level at finite places. Also the weights will be restricted to a cone inside the positive Weyl chamber. (This condition is parallel to the cone condition for geometric families and is important for similar reasons such as the uniform control of the analytic conductor). If we consider the level aspect then we fix a regular weight at infinity and consider a sequence of principal congruence subgroups of level $N\to \infty$.

The main theorem of~\cite{ST11cf} is a quantitative Plancherel equidistribution theorem for the local factors $\sigma_p$ of representations $\sigma\in\FmH$. Fix a test function $\varphi$ which is a Weyl invariant polynomial on the dual maximal torus of $H$.
For each prime $p$ large enough one can evaluate $\varphi$ against the unramified representations $\sigma_p$ of $H(\BmQ_p)$  and we have
\begin{equation}\label{plancherel}
	\sum_{\sigma \in \FmH(x)}
	 \varphi(\sigma_p)
	= \abs{\FmH(x)} \int  \varphi(\sigma) \mu^{\mathrm{pl}}_p(d\sigma)
	+ O(\abs{\FmH(x)}^{\delta} p^A)
\end{equation}
where $\mu^{\mathrm{pl}}_p$ is the unramified Plancherel measure on $\widehat{H(\BmQ_p)}$ and $\delta<1$.
The main term comes from the contribution of the identity on the geometric side of Arthur's cohomological trace formula~\cite{Art89}.
The remainder term comes from bounding the other orbital integrals. The multiplicative constant in $O(\cdot)$ is uniform in $p$ and $x$. This uniformity is a major difficulty in the proof because the number of conjugacy classes $\CmO$ to be considered on the geometric side is unbounded. In particular we have a weak control on the regularity of $\CmO$, it can for example ramify at several arbitrary large primes. We refer to~\cite{ST11cf}*{\S1.7} for a summary of the harmonic analysis techniques that we use to resolve this difficulty.

We deduce from~\eqref{plancherel} that each $\mu_p(\FmF)$ comes from the restriction of the Plancherel measure on $\widehat{H(\BmQ_p)}$. Precisely $\mu_p(\FmF)$ is the pushforward of $\mu^{\mathrm{pl}}_p$ under the functorial lift attached to $r:{}^LH \to \GL(n,\BmC)$. This is the assertion (i) of Conjecture~\ref{conj1}. The main term $t_\FmF(n)$ in the asymptotic~\eqref{quant} is expressed in terms of these $p$-adic densities.
We also get assertion (ii) and the global measure $\mu(\FmF)$ by inserting a more general test function $\varphi$ that is supported at finitely many places.

Maass forms are automorphic forms invariant under a maximal compact {subgroup} at infinity. They correspond to automorphic representations whose archimedean factors are spherical which is a condition that fits well in our formation of harmonic families. We expect the results to be similar to the case discussed above.
The classical case of Maass forms on $\GL(2)$ can be treated using the Selberg trace formula.
In higher rank the asymptotic Weyl's law is established in general by Lindenstrauss--Venkatesh~\cite{LV07}.  Weyl's law with remainder term and the quantitative equidistribution~\eqref{plancherel} are more difficult despite the harmonic analysis on the spherical unitary dual being well-understood~\cites{book:helg01,DKV}.\footnote{The difficulty is with the contribution of the continuous spectrum and in fact allowing noncongruence groups Weyl's law may fail~\cite{PS85}.} These difficulties revolve around the presence of Eisenstein series: notably there is not yet a satisfactory description of the residual spectrum for general groups.
The absolute convergence of the Arthur trace formula recently established by Finis--Lapid--M\"uller~\cite{FLM:spectral} is an important step forward. J.~Matz and the third-named author~\cite{Matz-Templier} have recently established the case of Maass forms on $\GL(n)$.


\subsection{Invariants of harmonic families}\label{sub:invariants}
We form the Sato--Tate measure $\mu_{\ST}(\FmF)=\lim\limits_{p\to \infty} \mu_p(\FmF)_{|T}$ in assertion (iii) of Conjecture~\ref{conj1}. Using the formula of Macdonald for the unramified Plancherel measure one can show this limit exists.\footnote{\label{f:nonsplit} This  holds literally if $H$ is a split group. For a general $H$ the Plancherel measure at a prime $p$ depends on the splitting behavior (it is \Lquote{Frobenian}). The \emph{average} of $\mu_p(\FmF)_{|T}$ over the primes $p< x$ as in~\eqref{STlimit} converges and assertion (iii) follows from Chebotarev equidistribution theorem.}
The measure $\mu_{\ST}(\FmF)$ coincides with the Sato-Tate measure attached to the image of ${}^LH$ viewed as a subgroup of $\GL(n,\BmC)$. This can be taken as the Sato-Tate group $H(\FmF)$ of the family, thus for harmonic families the existence of such a group is proven.

Next we examine the three indicators $i_1(\FmF)$, $i_2(\FmF)$ and $i_3(\FmF)$ in~\eqref{i}.  From now on we make the assumption that the representation $r:{}^LH\to \GL(n,\BmC)$ is irreducible which can be seen to be equivalent to $i_1(\FmF)=1$. Thus the family $\FmF$ is essentially cuspidal. This implies under the GRH that the functorial lift $r_* \sigma$ is cuspidal for most $\sigma\in \FmH$ which needs to be established by a separate unconditional argument. The strategy is to the relate the non-cuspidality of $r_* \sigma$ to the vanishing of certain periods of $\sigma$ (which is a well-studied and difficult problem, see the works of Jacquet, Jiang, Soudry and many others), that is that $\sigma$ is distinguished and then to show that this doesn't happen generically for almost all members $\sigma$ of $\FmF$. 

The indicator $i_2(\FmF)$ is either $1$ or $0$, depending on whether $r$ is self-dual or not. The indicator $i_3(\FmF)$ is denoted $s(r)$ in~\cite{ST11cf}.  It is the Frobenius--Schur indicator of $r$ which is either $-1$, $1$ or $0$, depending on whether $r$ is symplectic, orthogonal or not self-dual respectively. Thus the family $\FmF$ is essentially homogeneous if $r$ is irreducible and the homogeneity type is determined.


The rank $r(\FmF)$ is zero for harmonic families. This follows from the defining equation~\eqref{rank} and the Macdonald formula for the Plancherel measure which implies in every case the estimate $t_\FmF(p) = O(p^{-1})$, see~\cite{ST11cf}*{\S2}.
This vanishing of the rank reflects the fact that the central $L$-value (or the $L$-derivative if the root number is $-1$) is expected to vanish only for arithmetical reason which should happen only for a few exceptional members of the family $\FmF$. 

The root number is the most subtle of the invariants attached to the family $\FmF$. It is relevant for essentially symplectic families and corresponds to a decomposition
\[\mu_{\ST}(\FmF)=\mu_{\Sp,+}(\FmF)+\mu_{\Sp,-}(\FmF).\]
 For families in the level aspect the root number is related to the M\"obius function. See~\cite{ILS00} and the discussion below for the case of $\PGL(2)$.
In the weight aspect the root number could be dealt with along the lines of~\cite{ST11cf} although we have omitted the details there.

As we have noted repeatedly Conjecture~\ref{conj2} lies deeper.
Its formulation assumes the analytic continuation of the completed $L$-functions $\Lambda(s,\sigma,r)$ inside the critical strip in order to define the zeros. This is known in many cases notably via Rankin--Selberg integrals and the Langlands--Shahidi method. The functoriality conjecture of Langlands asserts that the $L$-functions should be attached to an isobaric representation $r_*\sigma \in \Aut(GL_n)$.

In this regard let us observe that under the Ramanujan conjecture for $\GL_n$ (resp. with a bound $\theta<\frac12$ towards Ramanujan), each of the local factors $L_v(s,\sigma,r)$ has no pole for $\MRe(s)>0$ (resp. $\MRe(s)> \theta$). Hence any zero $\rho$ with $\MRe(\rho)>0$ (resp. $\MRe(\rho)>\theta$) of the partial $L$-function $L^S(\rho,\sigma,r)=0$ cannot be cancelled by a potential pole of a local factor $L_v(s,\sigma,r)$ at $s=\rho$. The set of non-trivial zeros of $L^S(s,\sigma,r)$ (i.e., within the critical strip) will coincide with the set of zeros of $\Lambda(s,\sigma,r)$.
Thus Conjecture~\ref{conj2} only depends on the analytic continuation of the \emph{partial $L$-functions}. The formulation is robust because it is independent of the ramified factors $L_v(s,\sigma,r)$
(the analysis of which is the most delicate aspect in all known constructions of $L$-functions and the expected properties aren't fully established in many cases).

Once the above invariants $\mu_{\ST}(\FmF)$, $i_2(\FmF)$, $i_3(\FmF)$, $r(\FmF)$ and eventually $\mu_{\Sp,\pm}(\FmF)$ are found one can verify Conjecture~\ref{conj2} for a test function with restricted support. The size of the support depends directly on the quality of the estimate~\eqref{plancherel}.
The details are found in~\cite{ST11cf}*{\S12} while the Criterion 1.2 in~\cite{ST11cf} is the insight which has motivated our present formulation of Conjecture~\ref{conj2}.

We note that there is ample flexibility in choosing the spectral set $\FmH\subset \Aut(H)$. For example one can add harmonic analysis constraints at finitely many places. As soon as $\FmH$ is \Lquote{large enough}, the invariants of the family are independent of the choice and thus the Symmetry Type remains the same.
The analogue  for geometric families is to add congruences constraints on the parameters which is also very natural.

\subsection{Classical modular forms} \label{sub:gl2}
As mentioned above the case of $H=\PGL(2)$ is treated in~\cite{ILS00}.
One might wonder what an arbitrary parametrized spectral subset of $\Aut(H)$ should look like since our definition allows flexibility in choosing the local harmonic constraints.\footnote{In the context where $H$ is the unit group of a division algebra, P.~Nelson has recently proposed~\cite{Nelson:quaternion-statistics} conditions for certain test functions to isolate such `nice' spectral sets.} The problematic case of forms of weight $k=1$  is discussed in Section~\ref{sec:3}.
In this subsection we focus on the results of~\cite{ILS00} which  correspond to the spectral set of holomorphic cuspforms $S_k(N)$ of weight $k\ge 2$ and square-free level $N$ where either $k,N\to \infty$ with a possible additional average in dyadic intervals.

Suppose for simplicity that $r$ is the embedding $\SL(2,\BmC)\to \GL(2,\BmC)$ and denote by $\FmS$ the corresponding family of standard Hecke $L$-functions.
The conductor is $k^2 N$ and thus $\abs{\FmS(x)}$ which is the number of forms $f\in \FmS$ with $C(f)<x$ is asymptotic to $x$ up to a multiplicative constant.

Conjecture~\ref{conj1} holds for $\FmS$ as consequence of~\cite{ILS00} and the Plancherel equidistribution results~\cite{ST11cf} described in the previous subsections. The measure $\mu_{ST}(\FmS)$ is obtained from the conjugacy classes of $\SU(2)$ and hence coincides with the classical Sato-Tate measure.
If we let $T^1$ be the one-dimensional torus of $\SL(2,\BmC)$ and parametrize $T^1/W$ by
$\left(
\begin{smallmatrix}
e^{i\theta} & 0 \\
0 & e^{-i\theta}
\end{smallmatrix}
\right)$ with $0\le \theta \le \pi$, then
\[
\mu_{ST}(\FmS) = \frac{2}{\pi} \sin^2\theta \, d\theta.
\]

The indicators are given by $i_1(\FmS)=i_2(\FmS)=1$ and $i_3(\FmS) = -1$.
(More generally the Frobenius-Schur indicator of the $k$-th symmetric power representation $\SL(2,\BmC) \to \GL(k+1,\BmC)$ is equal to $(-1)^{k}$). Thus the family $\FmS$ is essentially symplectic and this is in accordance with the $SO(\infty)$ Symmetry Type.

To go further we decompose the family $\FmS = \FmS_+ \cup \FmS_-$ according to the root number being $+1$ or $-1$ respectively. The proportion of each piece is $\num{50}\%$.
The root number is $\varepsilon(f)=i^k \mu(N) \lambda_{f}(N)N^{\frac{1}{2}}$, so this statement is equivalent to cancellations in sums of the type $\sum\limits_{f\in \FmS(x)} \lambda_{f}(N)N^{\frac{1}{2}}$ which is an example of the M\"obius type sums discussed in~\eqref{mobius}. This sum can be analysed directly via the Petersson trace formula as in~\cite{ILS00} or alternatively using representation theory and the results in~\cite{ST11cf}. Above a prime $p\mid N$, the $p$-component of $f$ is tamely ramified with trivial central character and thus is either the Steinberg representation or a twist of the Steinberg representation by the unramified quadratic character; each representation carries $\num{50}\%$ of the mass of $\mu_p(\FmS)$ which comes from restriction  of the Plancherel measure  on $\widehat{\PGL_2(\BmQ_p)}$.

For $\Phi\in \CmS(\R)$ and $f\in \FmS$ we denote by $D(f,\Phi)$ the one-level distribution of the low-lying zeros of $\Lambda(s,f)$ (removing one zero at $s=\Mdemi$ if $f\in \FmS_-$). Then Conjecture~\ref{conj2} reads
\begin{equation}\label{ILS}
\frac{1}{\abs{\FmS_{\pm}(x)}}
\sum_{f\in \FmS_{\pm}(x)}
D(f,\Phi)
\rightarrow
\int_{-\infty}^{\infty}
\Phi(u) W_{\pm}^{(1)}(u)\, du,
\quad
\text{as $x\to \infty$}.
\end{equation}
In other words the Symmetry Type of $\FmS_+$ (resp. $\FmS_-$) is $\SO_{\even}(\infty)$ (resp. $\SO_{\odd}(\infty)$).

Unconditionally the asymptotic~\eqref{ILS} holds if the support of $\widehat \Phi$ is restricted to $(-1,1)$.
 Under the GRH for Dirichlet $L$-functions one can extend the support to $(-2,2)$.
This extension is significant because then the one-level density distinguishes between
the $\Sp(\infty)$, $\SO_{\even}(\infty)$ and $\SO_{\odd}(\infty)$ Symmetry Types since the distributions $W^{(1)}_+$ and $W^{(1)}_1$ agree in $u\in (-1,1)$ but split at $u=\pm 1$.


There are many interesting applications of $\GL(1)$ and $\GL(2)$ families, notably the non-vanishing of $L$-values, distribution of prime numbers, quantum chaos, subconvexity, equidistribution of arithmetic cycles and more. Here we have shown how to generalize the Symmetry Type with restricted support to higher rank families.
We view the low-lying zeros statistics as a first step towards these other
 arithmetic features and applications.

\subsection{GL(1) twists}\label{sub:twists}

We fix $\pi$ a cuspidal automorphic representation of $\GL(n)$ over $\Q$. If $\chi$ is a Dirichlet character we can consider the twist $\pi \otimes \chi$ which is again a cuspidal automorphic representation of $\GL(n)$.
In \S\ref{sub:n1} we have discussed $\GL(1)$ families, for example the family $\FmF^{(2)}$  of quadratic characters.  One can construct a parametric family
\[
\FmF=\set{\pi \otimes \chi,\ \chi\in \FmF^{(2)}}.
\]
As we have discussed in the remarks following the definition of families we allow one of the factor to be a singleton $\set{\pi}$ when considering the Rankin--Selberg product of families.

The quantitative equidistribution~\eqref{quant} is easily verified as well as the first two assertions of Conjecture~\ref{conj1}. The assertion (iii) however is as difficult as the individual Sato-Tate conjecture for $\pi$ itself.
We identify the $n$-dimensional torus $T$ with the diagonal of $\GL(n,\C)$ and thus with the product of $n$ copies of $\C^\times$. Assume the Sato-Tate conjecture holds for $\pi$ with a certain limit measure $\mu_{ST}(\pi)$ on $T$ and recall the Sato-Tate measure $\mu_{ST}(\FmF^{(2)})=\mu_B$ for $\FmF^{(2)}$ where $B=\set{1,-1}\subset \C^\times$.
We have a natural multiplication homomorphism $m:\C^\times \times T \to T$ given by pointwise multiplication of each coordinate. The assertion (iii) of Conjecture~\ref{conj1} holds and the Sato-Tate measure of the family $\FmF$ is the direct image
\begin{equation}\label{character-twist}
\mu_{ST}(\FmF) =
m_*(
\mu_B \times \mu_{ST}(\pi)
).
\end{equation}
Equivalently $\mu_{ST}(\FmF)$ is half the sum of $\mu_{ST}(\pi)$ and the image of $\mu_{ST}(\pi)$ under $t\mapsto -t$.
Note that since the family $\FmF$ is thin the average over the primes $p<x$ in~\eqref{STlimit} is critical (see also the Footnote~\ref{f:nonsplit} on page~\pageref{f:nonsplit} for another example).

Often $\mu_{ST}(\FmF)=\mu_{ST}(\pi)$, for example in the case that $\pi$ is a holomorphic modular form on $\GL(2)$ of weight at least two for which the individual Sato-Tate is known.
On the other hand the two measures may differ. The simplest example is when $\pi$ is a cubic Dirichlet character on $\GL(1)$ in which case $\mu_{ST}(\pi)$ is the Haar measure on the group $\set{1,e^{\frac{2i\pi}{3}},e^{\frac{4i\pi}{3}}}$ while $\mu_{ST}(\FF)$ is the Haar measure on $\set{
1,e^{\frac{i\pi}{3}},e^{\frac{2i\pi}{3}},-1,e^{\frac{4i\pi}{3}},e^{\frac{5i\pi}{3}}
}$.


In view of~\eqref{character-twist} and $\Mtr(-t)=-\Mtr(t)$ the indicators can be computed as
\[
\begin{aligned}
i_2(\FmF)&=
\int_{T} \Mtr(t)^2\, \mu_{ST}(\pi)(t)\\
 i_3(\FmF)&=
\int_{T} \Mtr(t^2)\, \mu_{ST}(\pi)(t)
\end{aligned}
\]
and $i_1(\FmF)=1$ since $\pi$ is cuspidal.
Thus we expect that $\FmF$ is essentially homogeneous and its homogeneous type is dictated by $\pi$. In fact we only need to know which of $L(s,\pi,\Msym^2)$ or $L(s,\pi,\wedge^2)$ has a pole at $s=1$, which is very little information about the Sato-Tate group $H_{\pi}$. So even if the Sato-Tate measure of $\pi$ remains mysterious we can verify the universality Conjecture~\ref{conj2} for $\FmF$ unconditionally, see~\cite{Rubin01}.
For example if $\pi$ is self-dual orthogonal, then $\FmF$ is essentially orthogonal and the Symmetry Type is $\Sp(\infty)$.

We can consider other $\GL(1)$ twists as for example the family $\FmF':=\set{\pi\otimes \chi}$ as $\chi$ ranges through all Dirichlet characters of conductor $q\le Q$ with $Q\to \infty$. Then the same analysis applies where we should replace $B$ by the full unit circle $S^1\subset \C^\times$. Thus we expect the Sato-Tate measure
\[
\mu_{ST}(\FmF') =
m_*(
\mu_{S^1} \times \mu_{ST}(\pi)
).
\]
The indicators are easier to compute in this case since we have $i_1(\FmF')=1$ and $i_2(\FmF')=i_3(\FmF')=0$. Thus the family $\FmF'$ is non self-dual and the Symmetry Type is $U(\infty)$ independently of any property of $\pi$. One simply uses  that $\pi$ is cuspidal and thus $L(s,\pi\times \widetilde \pi)$ has a simple pole at $s=1$ which controls all the restricted $n$-level densities universally. This is entirely analogous to the universality of high zeros found in~\cite{RS96}. This surprising universality and the behavior of the families $\FmF$ and $\FmF'$ fit nicely into our main conjectures.

We can analyse the previous example using the Sato-Tate group $H_\pi \subset \GL(n,\BmC)$, assuming it exists.
Then we would associate to the family $\FmF'$ the group $H(\FmF')$ generated by $H_\pi$ and $\BmC^\times$. In the same way that $\mu_{ST}(\pi)$ corresponds to $H_\pi$, we have that $\mu_{ST}(\FmF')$ corresponds to $H(\FmF')$.

Conversely we don't know what $H(\FmF')$ is unless we are willing to assume the existence of $H_\pi$.
 In fact
this example shows that if the family $\FmF$ is thin like this one, knowing $H(\FmF)$
is tantamount to knowing $H_\pi$ and so one may as well face having
to define $H_\pi$ conjecturally, for every $\pi$, if we want $H(\FmF)$ in general.

One expects that $H_\pi$ would be either a torus or semisimple. On the other hand $H(\FmF')$ obviously isn't and this immediately explains the vanishing of the indicators  $i_2(\FmF')=i_3(\FmF')=0$. In general a family whose Sato-Tate group has infinite center has to have an $U(\infty)$ Symmetry Type.


\subsection{Rankin--Selberg products}
In~\cite{DM06} Due\~{n}ez--Miller investigate an interesting example of a parametric family of $L$-functions obtained by a $\GL(2)\times \GL(3)$ Rankin--Selberg product.
Let $\pi$ be a fixed even unramified Hecke--Maass form on $\PGL(2)$. Consider the spectral set  of holomorphic cusp forms $f\in S_k(1)$ with $k\to \infty$. We can form the family
\[
\FmF := \set{\pi \times \Msym^2(f),\ f\in S_k(1)}
\]
which consists of $L$-functions of degree $6$. By the work of Kim and Shahidi functoriality is known in this case so $\FmF$ is a family of automorphic representations on $\GL(6)$. By construction all these forms are self-dual symplectic and the root number $\varepsilon(\Mdemi,\pi\times \Msym^2(f))$ can be verified to be $1$ for all $f$.

If we assume the Sato-Tate conjecture for $\pi$ then we can verify Conjecture~\ref{conj1} for $\FmF$. The measure $\mu_{ST}(\FmF)$ on the $6$-dimensional torus is associated to the subgroup $SU(2)\times PSU(2)$ of $U(6)$, where the embedding is given by $(\theta_1,\theta_2)\mapsto \theta_1\otimes \Msym^2 \theta_2$. Since
\[\Mtr(\theta_1 \otimes \Msym^2 \theta_2) = \Mtr(\theta_1)\Mtr(\Msym^2 \theta_2),
\]
 the indicators can be easily computed to be $i_1(\FmF)=i_2(\FmF)=1$ and $i_3(\FmF)=-1$.
 Thus the family is essentially symplectic as we expect.
In fact as usual we don't need to assume the full Sato-Tate conjecture for $\pi$ to compute these indicators, only the knowledge of the simple pole of $\Lambda(s,\pi \times \widetilde \pi)$ at $s=1$ suffices.

In~\cite{DM06} the $1$-level and $2$-level densities for a small restricted support are obtained unconditionally. This determines the Symmetry Type as $SO_{\even}(\infty)$ in Conjecture~\ref{conj2}.
This family $\FmF$ has the feature that each $L$-function has even functional equation without having to decompose a bigger family according to the root number, a feature which is present for any family with a $\Sp(\infty)$ Symmetry Type. Thus we can conclude following~\cite{DM06} that the Symmetry Type is not just a theory of signs of functional equations, which is also apparent in our Conjecture~\ref{conj2}.
More generally as studied in a subsequent paper~\cite{DM:convolving} the Symmetry Type has a certain predicted behavior under Rankin--Selberg product of families. This can also be explained by Conjecture~\ref{conj2} since if $\FmF_1$ and $\FmF_2$ are two essentially cuspidal homogeneous families we expect that $\FmF_1\times \FmF_2$ be homogeneous and in view of the properties of the Frobenius--Schur indicator that $i_3(\FmF_1\times \FmF_2) = i_3(\FmF_1)i_3(\FmF_2)$.

Another family that is constructed with Rankin--Selberg type integral consists of adjoint $L$-functions. For a family attached to the spectral set of Maass forms $\FmH$ on $\SL(3,\Z)$ this is studied recently by Goldfeld--Kontorovich~\cite{GK:low-lying} using their version of the Kuznetsov trace formula.
They consider the harmonic family $\FmF=(\FmH,\mathrm{Ad}_*)$ where $\mathrm{Ad}_*$ corresponds to the adjoint representation. The main result of~\cite{GK:low-lying} is that the family has Symmetry type $\Sp(\infty)$ when the density sums \eqref{F2} with $r=1$ are weighted by special values at 1 of $L$-functions of members of the family. (These weights are not expected to affect the Symmetry type.) This is consistent with our Conjecture~\ref{conj2} since $\FmF$ is a homogeneous family which is essentially orthogonal. Indeed if $\pi$ is a cuspidal automorphic representation on $\SL(3,\BmZ)$ then $L(s,\pi,\mathrm{Ad})$ is self-dual and orthogonal (it is always cuspidal because we are in full level thus $\pi$ is not a base change).

Actually this example generalizes nicely: let $H$ is any split connected quasi-simple group over $\BmQ$. Form the adjoint representation which is an $L$-map from ${}^LH$ to $GL_n$ where $n=\dim H$. Consider a generic spectral set $\FmH$ as above and the family $(\FmH,\mathrm{Ad}_*)$.
The adjoint representation is irreducible and it preserves the Killing form on $\TLie({}^LH)$ which is bilinear symmetric and non-degenerate.
Thus we expect almost all $L$-functions to be cuspidal and self-dual orthogonal thus the family to be essentially orthogonal.
Therefore according to Conjecture~\ref{conj2} we expect that any universal family of adjoint $L$-functions have Symmetry Type $\Sp(\infty)$.
For $H=\PGL(2)$ the adjoint representation is the same as the symmetric square and this is a result in~\cite{ILS00}.

The case of $H=\mathrm{PGSp}(4)$ is recently studied by Kowalski--Saha--Tsimerman~\cite{KST:Sp4}.
Namely they consider the spectral set $S_k^*(\Sp(4,\BmZ)) \subset \Aut(H)$ of Siegel cusp forms of weight $k\to \infty$.
Let $r$ be the degree four spin representation of ${}^LH = \Spin(5,\C)$.
We can form the family of $L$-functions
\[
\FmF := \set{L(s,F,r),\ F\in S_k^*(\Sp(4,\BmZ))}
\]
which by functoriality for classical groups are known  to correspond to automorphic representations of $\GL(4)$.

The main result of~\cite{KST:Sp4} is a (weighted) equidistribution result which is essentially related to Conjecture~\ref{conj1} for $\FmF$.
The measure $\mu_p(\FmF)$ is a (relative) Plancherel measure whose limit $\mu_{ST}(\FmF)$ exists as $p\to \infty$ and coincides with the Sato-Tate measure associated to the subgroup $r({}^{L}H)\subset\GL(4,\C)$.

One finds that $i_1(\FmF)=1$, thus the family is essentially cuspidal. The members $F\in S_k^*(\Sp(4,\BmZ))$ such that $L(s,F,r)$ is not cuspidal are precisely the Saito--Kurokawa lifts from $\SL(2,\BmZ)$. These form a (spectral) subset which is asymptotically negligible which confirms that almost all members of the family are cuspidal.

Next we have $i_2(\FmF)=1$ and $i_{3}(\FmF)=-1$ and thus the family is essentially symplectic.
In view of the isomorphism $\Spin(5,\BmC)\simeq \Sp(4,\BmC)$, the representation $r$ is self-dual symplectic which is consistent. The root number is $(-1)^k$ thus we expect according to Conjecture~\ref{conj2} an $\SO_{\even}(\infty)$ or $\SO_{\odd}(\infty)$ Symmetry Type, depending on the parity of the weight $k$.

The analysis of the low-lying zeros with a test function of restricted support is carried out in~\cite{KST:Sp4} but the results are altered by the presence of a weighting factor for each $F$. Since this weight is itself a \emph{central} value of $L$-function by a formula conjectured by B\"ocherer and Furusawa--Martin, it carries much fluctuation which apparently yields a symmetry which is not consistent with our conjectures.
 If these weights are removed we expect that this feature will disappear. Here this means that the weights which appear naturally from the application of the Petersson trace formula would need to be removed in order to interpret the symmetry type, see~\cite{Kowalski:families} for further discussions.

\subsection{Universal families}
For the universal family of all cuspidal automorphic forms on $\GL_n(\BmA)$ we expect that the Sato--Tate Conjecture~\ref{conj1} still holds. The measure $\mu_p(\FmF)$ is closely related to the Plancherel measure. Precisely for each integer $k\ge 0$, let $\mu^{\mathrm{pl}}[p^k]$ be the restriction of the Plancherel measure to the subset of representations in $\widehat{\GL_n(\BmQ_p)}$ of conductor $p^k$.  Then $\mu_p(\FmF)$ will be an explicit linear combination of the measures $\mu^{\mathrm{pl}}[p^k]$.

\begin{example}
		This can be verified for $n=1$, the universal family $\FmF$ of all Dirichlet characters, see also~\cite{Kowalski:families}. The total mass of $\mu^{\mathrm{pl}}[p^k]$ is $\varphi(p^k)$, the Euler function. A direct calculation shows that
\[
\mu_p(\FmF) = a \sum_{k=0}^\infty
\frac{1}{p^{2k}} \mu^{\mathrm{pl}}[p^k]
\]
 where $a=\frac{p^3}{(p-1)(p+1)^2}$.
\end{example}

Note however that for the ``family'' of forms of level $n!$ or product of consecutive primes $2\cdot 3 \cdot 5 \cdot 7...$, the Sato--Tate conjecture in the form~\eqref{quant} fails (as observed by Junehyuk Jung). The universality of the low-lying zeros in Conjecture~\ref{conj2} is still expected to hold here, but for deeper reasons. The case of families of Dirichlet characters can be verified directly, the case of $\GL(2)$ is done in~\cite{ILS00} and the general case is done in~\cite{ST11cf}.



\subsection{Deligne--Katz equidistribution and geometric families}\label{sub:katz}
In this subsection we consider geometric families. Our goal is to explain how to approach Conjecture~\ref{conj1} using monodromy groups. There are many technical issues that we ignore and we confine ourselves to an outline. 

We begin with a general geometric family as in the Definition in Section~\ref{section1}. Thus $\kW$ is an open dense subscheme of $\A_\Z^m$, and $f:X\to \kW$ is smooth and proper with integral fibers. To concentrate on examples of geometric nature, we assume the fibers to be geometrically connected. For any $w\in W:=\kW(\Z)\cap C$ we denote the fiber by $X_w$. This gives rise to a parametric family $\FmF$ of Hasse-Weil $L$-functions.

The local $L$-factor can be described using Grothendieck's $l$-adic monodromy theorem. (We need a result in $p$-adic Hodge theory when $p=l$ but it is harmless to assume $p>l$ for our purpose.) Let $\rho_w$ be the $\mathrm{Gal}(\overline \Q / \Q)$-representation acting on the space $H^d_{\acute{e}t}(X_w\times_\Q \ol{\Q},\Q_l)$. For any prime $p$ we consider the Weil-Deligne representation
\[
r_{w,p} := \iota \WD_v(\rho|_{\Gal(\ol{\Q}_p/\Q_p)}),
\]
see Appendix~\ref{sub:Hasse-Weil} for details. Also let $\pi_{w,p}:=\rec^{-1}(\rho_{w,p})\otimes|\det|^{d/2}$ viewed as an element of $\widehat{G(\Q_p)}$ where $G=\GL_n$. (As remarked in the previous section, the fact that $\pi_{w,p}$ is unitary is conditional on the weight-monodromy conjecture if $X_w$ has bad reduction.) The local $L$-factor at $p$ is given by
\[
L(s,\pi_{w,p})=L(s,r_{w,p})
= \det (1-\Frob_p p^{-s}|V^{I_p}\cap \ker N)^{-1}
\]
where $V$ is the underlying space of $r_{w,p}$ and $N$ is the corresponding nilpotent operator.

As a preliminary step we examine the ramification of the representations $\pi_{w,p}$. If $\pi_{w,p}$ is ramified then $p$ is a prime of bad reduction for $X_w$ and also $D(w)\equiv 0 \pmod{p}$, where $D$ is the discriminant function of the family. Conjecture~\ref{conj1} is rather precise because the assertions (i) and (ii) include the statistics of the ramified representations. The depth of the representations $\pi_{w,p}\in \widehat{G(\Q_p)}$ is bounded by a constant~\cite{ST:field-of-rat}*{\S3} independent of $w,p$ because its field of rationality is $\Q$. 

For each unramified  $\pi_{w,p}$ we obtain an element $t_{w,p}\in T/W$. A crucial observation is that it depends only on $w$ modulo $p$. Thus the measure $\mu_p(\FmF)|_T$ (and more generally $\mu_p(\FmF)$) is atomic, in fact supported on a finite subset of $T/ W$. It is given explicitly by the following sum of Dirac measures:
\begin{equation}\label{mup=sum}
\mu_p(\FmF)|_T = \frac{1}{|\kW(\F_p)|}\sum_{\substack{w\in \kW(\F_p),\\ D(w)\not\equiv 0}} \delta_{t_{w,p}},
\end{equation}
where the sum has been restricted to those $w$ such that $\pi_{w,p}$ is unramified by demanding that $D(w) \not\equiv 0 \pmod{p}$. It implies by the Lang-Weil bound~\cite{Lang-Weil},
\begin{equation}\label{mupF}
		\mu_p(\FmF)(T) = 1 - O\left(\frac{1}{p}\right).
\end{equation}

In view of \eqref{mupF} the ramified representations play no role in the assertions (iii) and (iv) of Conjecture~\ref{conj1} and hence also in the construction of $\mu_{\textrm{ST}}(\FmF)$ which is our main interest. Thus from now on we shall focus on~\eqref{mup=sum} and those representations $\pi_{w,p}$ which are unramified.

The analysis involves sets of integer points $w=(w_1,\ldots,w_m)$ in sectors $W$ in $\Z^m$ in regions defined by a homogeneous polynomial which approximates the conductor, for example a height condition that $w$ lies in a large box (that is each $w_i$ lies in an interval). The sectors defining $C$ are chosen to make these sets finite by avoiding the projective zero locus of the discriminant $D$. The assertion~\eqref{mup=sum} is deduced from the convergence:
\[
		\frac{1}{|\FmF(x)|}\sum_{\substack{w\in \kW(\Z)\cap C,\\ |w_i|^{d_i}<x, \forall i}}  \delta_{t_{w,p}} \rightharpoonup
	\frac{1}{|\kW(\F_p)|}\sum_{\substack{w\in \kW(\F_p),\\ D(w)\not\equiv 0}} \delta_{t_{w,p}}	,
\]
which follows from the fact that $t_{w,p}$ depends only on $w$ modulo $p$.

The above reasoning is the key arithmetic input. And indeed this argument occurs often in number theory such as in the circle method.
This localization away from the zero locus of $D$ makes the problem easier and in general it forces us to count the parametrized elements $\pi_w$ in the family with some natural multiplicity.\footnote{We note that we make analogous simplifying assumption in the case of harmonic families, see~\S\ref{sub:harmonic}, where we have allowed some mild weights such as $\dim(\pi_v)^{U_v}$ which doesn't change the final answer but makes the problem easier to analyse with the trace formula.}

To establish assertions (iii) and (iv) of Conjecture~\ref{conj1} it remains to study the measures $\mu_p(\FmF)|_T$ and thus we are reduced to a problem over finite fields. The reduction is possible because we have chosen $\kW$ to be affine in the definition of geometric families. In fact we see from the argument that we could relax this assumption somewhat, but not entirely see~\S\ref{sub:limitation} below.

It is convenient to formulate the problem over finite fields by introducing the sheaf $\kG:=R^df_* \ol{\Q_l}$, which is the \Lquote{$H^d$ along the fibers $X_w$}. It is a lisse $\ell$-adic sheaf over $\kW$ of rank $n$.
The Grothendieck base change theorem implies that there is an action of the arithmetic fundamental group $\pi_1(\kW)$ on a finite dimensional $\ol{\Q_l}$-vector space which can be identified with the cohomology group of the fibers~\cite{book:KS}. Specifically there is a linear action by automorphism which yields the monodromy representation $\pi_1(\kW)\to \GL(n,\overline{\Q_l})$, which is well-defined up to conjugation. The geometric fundamental group $\pi_1^{\geom}(\kW)$ is a normal subgroup of $\pi_1(\kW)$, and we denote by $G_{\geom}$ the Zariski closure of its image. By a theorem of Deligne $G_{\geom}$ is semisimple. The Zariski closure of the image of the arithmetic fundamental group $\pi_1(\kW)$ is denoted $G_{\arith}$.\footnote{Here we are assuming as in~\cite{Katz:elaboration} geometric connectedness.
}
Thus $G_{\geom}$ is a normal subgroup of $G_{\arith}$ and from now on we impose the hypothesis that $G_{\arith}\subset \mathbb{G}_mG_{\geom}$, see the recent article of Katz ~~\cite{Katz:elaboration} for details, cf. Hypothesis (H) there. This essentially amounts to a purity assumption on the sheaf $\kG$, which gives a uniform control on the size of Frobenius eigenvalues.


For each prime $p$ and $w\in \kW(\F_p)$, the image of Frobenius under the monodromy representation lies in $G_{\arith}$. Thanks to the hypothesis above, we can rescale it by a scalar and obtain an element $\Frob_{w,p} \in G_{\geom}$ well-defined up to the choice of an $l$-adic unit and up to conjugation. Moreover by purity all the eigenvalues of $\iota\Frob_{w,p}$ lie on the unit circle and therefore $\iota\Frob_{w,p}$ may be viewed up to conjugation as an element of $B_c$, the maximal compact subgroup of $G_{\geom}$, again we refer to~\cite{Katz:elaboration} for details.

We form the probability measure
\[
\mu_p(\kG):= \frac{1}{|\kW(\F_p)|}\sum_{w\in \kW(\F_p)} \delta_{\Frob_{w,p}}
\]
on $B_c^\#$.
The key point of these constructions is that the pushforward of $\mu_p(\kG)$ under $B_c^{\#} \to T_c/W$ coincides up to $O\left(\frac1p\right)$ with the measure $\mu_p(\FmF)|_T$ defined in~\eqref{mup=sum} above.

It remains to let the prime $p\to \infty$. The equidistribution of the measures $\mu_p(\kG)$, with respect to the Haar measure of $B$, is Katz's variant of Deligne's equidistribution theorem, see \cite{Katz:elaboration} and~\cite{book:KS}*{\S9}. It is important here that it can be proven that the monodromy depends only on the topology of the family $X\to \kW$. In other words the geometric fundamental group is independent of $p$ for $p$ large, see~\cite{Katz:elaboration}*{Thm.~2.1}.

Specifically we apply Theorem~5.1 of~\cite{Katz:elaboration} (with all $n_i$ equal to 1) to the sheaf $\kG$, which is $\iota$-pure by \cite[9.1.15]{book:KS}, to obtain that
\begin{equation}\label{katz}
\mu_p(\kG) \rightharpoonup \mu_{ST}(\FmF),\quad
\text{as $p\to \infty$,}
\end{equation}
where $\mu_{ST}(\FmF)$ is the pushforward of the Haar measure under $B_c^\# \to T_c/ W$. Note that the base scheme $S$ for us is of the form $\mathrm{Spec} \Z[1/N]$ and therefore the Hypothesis (AFG) in~\cite{Katz:elaboration} involves removing finitely many primes $p$. This finishes the outline of the proof of the assertions (iii) and (iv) of Conjecture~\ref{conj1} for $\FmF$.

For example for the family of all elliptic curves which we have discussed in~\S\ref{sub:elliptic}, the equidistribution theorem is an early result of Birch~\cite{Birch:how}. The example of $1$-parameter families of hyperelliptic curves of genus $g$ is treated in~\cite{book:KS}, where we have $G_{\geom}=\Sp(2g)$ and $G_{\arith}=\GSp(2g)$. Another interesting example is the universal family of smooth projective hypersurfaces of given dimension and degree, which is also in~\cite{book:KS}. Finally the above equidistribution applies to the Dwork families discussed in~\S\ref{sub:dwork}.

Conjecture~\ref{conj2} can be established for $\FmF$ for test functions of limited support and conditionally on the modularity conjecture for the $X_w$. Both for harmonic families (see \S\ref{sub:invariants}) and for geometric families we have attached a group $H(\FmF)$ such that the associated Sato-Tate measure $\mu_{ST}(\FmF)$ is computed in terms of $H(\FmF)$. As we observed earlier the measure $\mu_{ST}(\FmF)$ need not determine the group $H(\FmF)$ uniquely, however there is a natural choice which comes from the method of proof of Conjecture~\ref{conj1}, namely $H(\FmF):=r({}^LH)$ for harmonic families and $H(\FmF):=G_{\geom}$ for geometric families.

Serre has recently put forward a Sato-Tate conjecture for schemes which is related to the above discussion. 
Let $X\to \kW$ be a scheme of finite type. If $\kW$ is a point then this is the usual Sato-Tate conjecture for the Hasse--Weil $L$-functions attached to $X$. If $\kW$ satisfies some suitable conditions it is a direct consequence of~\eqref{katz} as explained in~\cite{Katz:elaboration} because it asks for the convergence for $x\to \infty$ of the average for $p^r<x$ of the measures $\mu_{p^r}(\kG)$. There are differences of this to our Sato-Tate conjecture for families: one being that the Sato-Tate conjecture for scheme is expected to be true for any base $\kW$ (and is proven in~\cite{Katz:elaboration} under mild assumptions if $\kW$ is not a point), whereas it is easy to construct counterexamples to our Conjectures~\ref{conj1} and \ref{conj2} for families if the base $\kW$ were arbitrary (see~\S\ref{sub:limitation}).

\subsection{Prospects}\label{sub:prospects} 
Under certain assumptions we have verified for the above families the concepts introduced in Section~\ref{section1}. 
It is desirable to lift these assumptions as much as possible since this would strenghten our knowledge and make certain results unconditional. 
We summarize here the nature of these issues and give some plausible outlook of how some could be addressed in future work. 
We shall focus solely on the Sato-Tate equidistribution for families as formulated in Conjecture~\ref{conj1}.  

For general harmonic families, the Sato-Tate equidistribution for families implies working with general test functions, which raises important questions on the global harmonic analysis of the trace formula.
One such question is formulated in~\cite{FLM:mult} in the context of limit multiplicities and concerns a uniform estimate on the winding number of normalizing scalars of intertwinning operators.
Another challenge concerns the description of the residual spectrum which is known for $\GL(n)$ and used crucially in establishing quantitative error terms in the Weyl's law~\cite{LM09,Matz-Templier}.
These and related problems now seem within reach in the context of classical groups from the work of Arthur and others.

Local harmonic analysis and representation theory of $p$-adic groups and real Lie groups also play a major role in Conjecture~\ref{conj1}.
One would like to capture a portion of the spectrum that is as fine as possible. 
Over the reals this means discrete series versus stable packets and short spectral windows for Maass forms.
For $p$-adic groups this means working with congruence subgroups beyond principal towers, see e.g.~\cite{Finis-Lapid:congruenceII}, and possibly working with a single supercuspidal representation, a question discussed in~\cite{t:trace-character} which will appear in this proceedings volume. 
Another property concerns uniform control on the matrix coefficients of intertwinning operators, which is studied in~\cite{MS:gln} over the reals and in~\cite{FLM:degrees} over the $p$-adics.
Finally the analytic conductor of representations, which is used in the present formulation of Conjecture~\ref{conj1}, is difficult to define in complete generality. 
For this it is essential to clarify the relation between depth and conductor, see~\cite{Kala:phd} for work in this direction, and it would be important to improve our understanding of the local Langlands correspondence in the tame case.

For geometric families it is a difficult problem in each specific example to identify the monodromy group. 
Also it is difficult to make the parametrization $F$ one-to-one; this is related to the implementation of the square-free sieve, which a major step in the work of Bhargava on counting number fields with bounded discriminant.
Analogously to the question of depth versus conductor mentioned above for automorphic representations, there is a question of the relation between height and conductor for Hasse-Weil $L$-functions.

\section{Section 3}\label{sec:3}

		In this section we give some ``families'' of automorphic forms that do not fit into our prescription in Section~\ref{section1}. While some of these are natural and Conjectures~\ref{conj1} and \ref{conj2} probably apply to them, they lack parametrizations and hence any known means of study and hence remain very speculative.





\subsection{Limitations}\label{sub:limitation} We begin by pointing to limitations in forming families.
The base space $\kW$ of parameters in our definition of a geometric family is allowed to be $\mathbb{P}^m/\Q$, $\A^m/\Z$ or products of such. Unlike the algebro-geometric setting of families over finite fields, we cannot allow a general base $\kW$ which is defined by equations over $\Z$ (or $\Q$). According to the solution of Hilbert's 10th problem~\cite{Matiyasevich:hilbert} one cannot say much about such sets $\kW(\Z)$, for example deciding if they are finite or not, and in general these sets may be unwieldy (see the example below). In particular the averages \eqref{quant}, or for that matter any other statistics associated with the family, need not exist. What would suffice for $\kW(\Z)$ in order for us to analyze the family to the extent that is described in Sections \ref{section1} and \ref{sec:examples} is that $\kW$ be ``strongly Hardy--Littlewood'' in the sense  of~\cite{Borovoi-Rudnick}.

The same difficulty arises if we try to perform simple Boolean operations on our families. If $\FmF_1=(W_1,F_1)$ and $\FmF_2=(W_2,F_2)$ are two parametric families in $\Aut(G)$ then a natural parametric definition of their intersection is $\FmF_{12}=(W_{12},F_{12})$ where $W_{12}=\{(w_1,w_2):F_1(w_1)=F_2(w_2)\}\subset W_1\times W_2$ and $F_{12}((w_1,w_2))=F_1(w_1)$ ($=F_2(w_2)$) for $(w_1,w_2)\in W_{12}$. Note that if $F_1$ and $F_2$ are embeddings (so that $F_1(w_1)$ an $F_2(w_2)$ are parametrized sets in $\Aut(G)$) then $\FmF_{12}$ parametrizes $\FmF_1(W_1)\cap \FmF_2(W_2)$. The problem is that $W_{12}\subset W_1\times W_2$ encodes a general diophantine set and again we are dealing with unwieldy sets for which the various statistical averages over the family need not exist.

  A concrete example of the above where we allow various operations on a parametric family is the following: Let $R\in \N$ be a recursive set~\cite{Matiyasevich:hilbert}. There is a polynomial $P=P_R\in \Z[W_1,...,W_{10}]$ such that $P(\Z^{10})\cap \N=R$ (see~\cite{Matiyasevich:hilbert}). Consider the parametric family $\FmF$ in $\Aut(\GL_1)$ given by
   $$\FmF:~X^2=p(W_1,...,W_{10})$$
   so that
   $$F((w_1,...,w_{10}))=\chi_{D(w_1,...,w_{10})}$$
   where $D(w_1,...,w_{10})$ is the square-free part of $p(w_1,...,w_{10})$ and $\chi$ the Dirichlet character corresponding to the quadratic field $\Q(\sqrt{p(w)})$. Then $\FmF=(W,F)$ is a parametric family in our sense and the discussion in Sections \ref{section1} and \ref{sec:examples} applies to it. However if we consider the image $T=F(\Z^{10})$ in $\Aut(\GL_1)$ and impose the condition that the field corresponding $F(w)$ is real (that is we intersect $T$ with $\N$) then we arrive at the subset $R$ of $\N$, realized as a subsect of $\FmF^{(2)}$. The set of recursive subsets of $\N$ is very general and certainly any statement such as \eqref{quant} will not hold for such a general $R$ (when ordered by height).

\subsection{Fields of rationality}\label{sub:rationality} In this section we introduce a construction of families via field of rationality. Let $\pi$ be an automorphic representation of $\GL_n(\A)$. The field of rationality $\Q(\pi)$ for $\pi$ is by definition the fixed field in $\C$ under $$\{\sigma\in \mathrm{Aut}(\C): \pi^{\sigma}\simeq \pi\}$$ where $\pi^{\sigma}:=\pi\otimes_{\C,\sigma} \C$. A well known conjecture states that $[\Q(\pi):\Q]<\infty$ if and only if $\pi$ is algebraic in the sense of Clozel \cite{Clo90}. (These notions and the conjecture extend to arbitrary connected reductive groups, cf. \cite{BG}.)

Let $\mathfrak{F}=(\mathfrak{H},F)$ be a harmonic family as in Section \ref{section1}. For a number field $K$ (as a subfield of $\C$) define $\mathfrak{F}_{\subseteq K}$ to be the subset consisting of $\pi\in\mathfrak{F}$ such that $\Q(\pi)\subseteq K$. Similarly for an integer $A\ge 1$ define $\mathfrak{F}_{\le A}:=\{\pi\in\mathfrak{F}: [\Q(\pi):\Q]\le A\}$. Observe that each of $\mathfrak{F}_{\subseteq K}$ and $\mathfrak{F}_{\le A}$ is supposed to contain only algebraic members by the conjecture just mentioned. If $\mathfrak{F}$ is ramified at only finitely many primes then $\mathfrak{F}_{\subseteq K}$ and $\mathfrak{F}_{\le A}$ are conjectured to be finite sets, cf. \cite[Conj 5.10]{ST:field-of-rat}, and verified to be finite when  $G$ is a general linear group or a quasi-split classical group. (See Theorem 1.6 and Corollary 6.8 of \cite{ST:field-of-rat}.)

\begin{example}
  In the setup for harmonic families take $H=G=\GL_2$. Let $\mathfrak{H}$ be the family of all cuspidal automorphic representations $\pi$ of $\GL_2(\A)$ such that $\pi_\infty$ is the discrete series of lowest weight (so that $\pi$ correspond to classical modular forms of weight 2). Suppose that $F$ comes from the identity $L$-morphism $r$. Then $\mathfrak{F}_{\subseteq \Q}=\mathfrak{F}_{\le 1}$ is nothing but the family of all normalized cuspforms of weight 2 whose Fourier coefficients are rational numbers.
\end{example}

The family $\mathfrak{F}_{\subseteq \Q}$ in the example is identified with the family of all elliptic curves over $\Q$, cf. Appendix \ref{sub:Hasse-Weil} below. The family corresponds to the moduli stack of elliptic curves over $\Q$ or a moduli scheme if a suitable level structure is added. So this example almost fits in the framework of geometric families considered earlier, to which the two main conjectures apply.
  This leads us to the question as to when the families $\mathfrak{F}_{\subseteq K}$ and $\mathfrak{F}_{\le A}$ can be realized as geometric families. Moreover we may ask

\begin{question}
  Suppose that the family $\mathfrak{F}_{\subseteq K}$ (resp. $\mathfrak{F}_{\le A}$) has infinite cardinality. Are Conjectures \ref{conj1} and \ref{conj2} true for the family $\mathfrak{F}_{\subseteq K}$ (resp. $\mathfrak{F}_{\le A}$)?
\end{question}

  To shed light on the question, let us pursue the connection with geometric families further when the family $\mathfrak{F}_{\subseteq \Q}$ is constructed as in the above example except that the weight is a general integer $k\ge 2$, following \cite{PR:Calabi-Yau}. (Also see \cite[\S7.2]{Kha10}.) A conjecture of Paranjape and Ramakrishnan states that each $\pi\in \mathfrak{F}_{\subseteq \Q}$ should be associated with a two-dimensional $\Gal(\ol{\Q}/\Q)$-subrepresentation of $H^{k-1}(X_\pi\times_\Q \ol{\Q},\ol{\Q}_l)$ for some Calabi-Yau variety $X_\pi$ over $\Q$ of dimension $k-1$ (such that the two-dimensional piece should be cut out by the part with Hodge numbers $(k-1,0)$ and $(0,k-1)$). If true, this suggests that $\mathfrak{F}_{\subseteq \Q}$ might be a family of 2-dimensional motives appearing in the family of $H^{k-1}$-cohomology arising from a family of $(k-1)$-dimensional Calabi-Yau varieties. When $k=2$ this reduces to the discussion of the family of elliptic curves over $\Q$ above. In case $k=3$, where all $\pi$ are of CM type and $X_\pi$ are K3 surfaces, see \cite{ES:weight3} for a recent result due to Elkies and Sch\"utt. A partial result towards the general case is worked out in \cite{PR:Calabi-Yau}. However it is known that there are only finitely many $\pi$ which are of CM type, correspond to a weight 3 cuspform, and have $\Q$ as field of rationality, and similarly for all odd $k\ge 3$ under the GRH, cf. \cite[\S3]{ES:weight3} for more details.  So the assumption of the above question is not superfluous. In fact the authors do not know a criterion for $\mathfrak{F}_{\subseteq K}$ to be infinite. 

  %

More generally these conjectures about other rationality for algebraic representations all point to geometric families again.  So philosophically perhaps many families obtained by specifying the field of definition are already included in our geometric families. (However it may be too bold to predict that all such families obtained by constraining the field of rationality can be constructed via geometry. For instance the case of $\GL(n)$ for $n\ge 3$ is unclear.) On the other hand, we note that a result on the degree of the field of rationality by two of us (\cite{ST:field-of-rat}) can be interpreted as the following statement: a harmonic family cannot be defined by a geometric construction, at least when the components at infinity are discrete series, because then the degree of the field of rationality would be bounded.

There are other examples such as the family of all Maass forms of eigenvalue $\frac{1}{4}$, say with integer coefficients. A letter~\cite{Sarnak:Maass-integer}, extended in~\cite{Brumley:Maass-integer} shows that these  forms are the same as certain  Galois representations with a given $H$-type (see below). So this family too can be thought of in two ways.

\subsection{Local conditions with measure zero}\label{sub:measure-zero}

In the construction of harmonic families we allowed ourselves to restrict a local component $\pi_v$ to a nice subset $B_v\subset \widehat{H(\Q_v)}$ only for $B_v$ of \emph{positive} Plancherel volume. It is of interest to study some cases where $B_v$ has measure zero. In doing so our main tool for studying the family, namely the trace formula, cannot be used effectively to isolate members of the family.

An important special case is to take $\pi_\infty$ in a specified finite subset. For a fixed irreducible algebraic representation $\xi$ of $H$ over $\C$, take $B_\infty$ to be the set of $\pi_\infty\in \widehat{H(\R)}$ such that $\pi_\infty$ is cohomological for $\xi$, namely $\pi_\infty\otimes \xi$ has nonzero Lie algebra cohomology in some degree. Then $B_\infty$ is a finite set and often has Plancherel measure zero, for instance when $H=GL_n$ for $n\ge 3$ and $\xi$ is arbitrary. Then $\pi\in \Aut(H)$ is such that $\pi_\infty\in B_\infty$ captures the information about the cohomology of the corresponding locally symmetric space for $H$ with coefficients in a local system arising from $\xi$.  One could refine the above choice of $B_\infty$ by taking $B_\infty$ to be a singleton $\pi_\infty\in \widehat{H(\R)}$ which is cohomological for $\xi$ but not a discrete series.
As a further generalization of the special case above, one can take $B_\infty$ to be a finite set consisting of $\pi_\infty\in \widehat{H(\R)}$ which are C-algebraic in the sense of \cite[Def 2.3.3]{BG}. Roughly speaking, it means that the infinitesimal character of $\pi_\infty$ is integral after twisting by the half sum of all positive roots of $H$. For example we get the family of all weight 1 cuspforms and the family of all Maass cuspforms with Laplace eigenvalue 1/4 when $H=\GL(2)$ and $B_\infty$ is a suitably chosen singleton.

In all of these cases it is already a difficult  problem to enumerate $\mathfrak{H}$ as analytic conductor grows, in other words to study the asymptotic growth of \eqref{weyllaw}. The answer to the last is well known when $B_\infty$ consists of discrete series (and thus has a positive volume) by work of de George-Wallach~\cite{DeGeor-Wall}.\footnote{A general \emph{uniform} such limit multiplicity theorem has been derived recently in~\cite{BG7:betti11}.} In the case at hand concerning these families for which $B_\infty$ is as above, there have been some conjectures and results concerning the sizes of these sets, see~\cite{SX91}, \cite{CE09}, \cite{Mar12}. Here we take a step further to pose the question of whether our main conjectures (Conjectures \ref{conj1} and \ref{conj2}) are true for such families. The same question can be asked when we prescribe constraints at finite places by subsets of Plancherel measure zero.

\subsection{Universal $H$-types} As discussed above any of our pure families $\FmF$ has an $H$-type associated with it, namely an $H$ such that $\mu_{\ST}(\FmF)=\mu_H$. Conversely one might try form universal families with a given $H$-type. Given $H$, the set of $\pi$'s in $\Aut_\cusp(G)$ with $H_\pi=H$ would be such a family, or we could impose this condition on $\pi$'s in any one of our families. There are some basic difficulties with such a construction. The first is that we don't know how to define $H_\pi$ in general. To begin with we can get around this problem by restricting to $\pi$'s which are algebraic. The second problem is more serious and this is, in any generality we have no means of understanding such an $\FmF_H$ and even the simplest requisite~\eqref{weyllaw} is mysterious. Nevertheless it would seem safe to expect that the $H$-type of $\FmF_H$ is $H$, and that Conjectures~\ref{conj1} and~\ref{conj2} would hold for any rich enough such $\FmF_H$ (for example it should at least be an infinite set). A numerical study of such ``families'', even for $\GL_2$-forms, would be revealing. The difficulty with a theoretical study of such $\pi$'s is closely related to (but easier since we only ask asymptotic questions) to the analytic problem of recognizing $\pi$'s in $\Aut_\cusp(G)$ with a given $H_\pi$ that is raised by Langlands in his ``Beyond Endoscopy'' paper~\cite{Lang:beyond}.



While we can't attack these $H$-type families, we can in all cases (at least where the Noether conjecture is known) produce geometric parametric subfamilies of any of these types. In many cases these subfamilies are probably close to being a positive proportion of the $H$-types. In fact one of the standard approaches to the inverse Galois problem for special finite $H$'s is to make an $H$-extension of $\Q(T_1,T_2, \ldots ,T_m)$ and then to specialize the $t$'s and use Hilbert irreducibility by counting (see~\cite{Serre:Mordell-Weil}). This very construction is a geometric parametric family according to our definition and of course it gives a large subfamily of such an $H$-type in our context.

There are some $H$'s for which $\FmF_H$ can be studied, primarily using class field theory. For $G=\GL(1)$ and $H$ a finite cyclic subgroup of $\C^\times$, $\FmF_H$ consists of all Dirichlet characters of order $|H|$ (for $|H|=2$ this is the family $\FmF^{(2)}$ from \S\ref{sub:n1}). Conjecture~\ref{conj1} is established without much trouble and $\mu_\ST(\FmF_H)=\mu_H$ and for $|H|\ge 3$, $i_2(\FmF_H)=0$ and the symmetry type is $U(\infty)$. Conjecture~\ref{conj2} has been established for test functions of restricted support and numerically for $|H|=3$ (\cite{GZ:one-level,DFK:vanishingunitary}).

For $G=\GL(2)$ an interesting family related to $H$-types, with $H$ not fixed but varying itself over a class of groups, was constructed by Hecke. Namely $\pi$'s which are holomorphic cusp forms of weight $1$ for which $H_\pi$ is (finite) dihedral. One can study a refined version of Conjectures~\ref{conj1} and \ref{conj2} for this family by collecting these forms into smaller packets which correspond to Hecke characters of the class group of $\Q(\sqrt{D})$, where $D\to -\infty$. This was done in~\cite{FI03} who show that the symmetry type is $\Sp(\infty)$. From our point of view this is ``clear'' since $H_{\FmF}$ is a dihedral subgroup of $\GL_2(\C)$ and in particular has Frobenius--Schur indicator equal to $1$. Other than using class field theory and specifically $1$-dimensional characters, we know of few examples where universal families of $H$-types can be studied.

\subsection{Closing comments} \label{sub:closing}
There are obvious variations on these constructions. We can combine number field (geometric) families and harmonic families. For example let $\{K_i\}_{i\in I}$ be a family of number fields over $\Q$ of fixed degree $d$ such that $\disc(K_i)\rightarrow \infty$. A further option is to require that in addition that $K_i$'s have isomorphic Galois groups, that they satisfy a constraint on primes of ramification, or some other reasonable properties. Let $H$ be a connected reductive group over $\Q$, with an $L$-group representation $r: {}^L H \rightarrow GL(m,\C)$. The latter gives rise to an $L$-group representation $R: {}^L (\Res_{K_i/\Q} H) \rightarrow GL(md,\C)$ by applying $r$ on each copy of the dual group of $H$. The functorial lift corresponding to $R$ is the functorial lift with respect to $r$ over $K_i$ followed by the automorphic induction from $GL_m$ over $K_i$ to $GL_{md}$ over $\Q$. The resulting family $\FmF$ is a family of automorphic $L$-functions of degree $md$. If functoriality for $r$ (over each $K_i$) is known then we may think of $\FmF$ as a family of automorphic representations of $\GL(md,\A)$ whose standard $L$-functions are as above. Sometimes it happens that every $L(s,\pi,R)$ factorizes as a product of $L$-functions and has a certain factor in common. In that case we may as well remove the common factor altogether. This construction yields examples which are not covered by families of the first chapter.


Finally note that for any of our parametric families one can impose further restrictions in exhausting $\FmF$ or placing arithmetic conditions on the conductors. For example one can collect the $\pi$'s in $\FmF$ in shells of given conductor (going to infinity) if these sets are large, or one can restrict to $\pi$'s in $\FmF$ with conductor a prime number. We view these as simple variations of our formation of families, albeit often technically more problematic. We have emphasized families which are cuspidal and pure, however mixed types arise naturally enough. A good example is that of Dedekind zeta functions of quartic field extensions of $\Q$. For these a positive proportion have Galois closure $S_4$ (as in~\S\ref{sub:number-fields}) but there is also a positive proportion with Galois group $D_4$ whose invariants are quite different (see~\cite{Bha05}).
\appendix

\section{Hasse--Weil $L$-functions}\label{sub:Hasse-Weil}

  Here we recall the definition of the Hasse-Weil $L$-function \eqref{e:Hasse-Weil} and the modularity conjecture. The modularity conjecture (Conjecture \ref{conj:modularity} below) states that the $L$-functions arising from algebraic varieties over $\Q$ should be automorphic $L$-functions. In fact we will
explain how $L$-functions are attached to $l$-adic Galois representations, in particular the \'etale cohomology space appearing in \eqref{e:Hasse-Weil}. To do so we recall the local Langlands correspondence for general linear groups in order to be precise about the matching of $L$-functions at ramified places. We also reformulate the modularity conjecture as a bijective correspondence between certain $l$-adic Galois representations and automorphic representations preserving $L$-functions, incorporating observations by Clozel and Fontaine-Mazur.
  The reader is referred to \cite{Tay04} for an excellent survey of many topics discussed in this appendix.

  Let $p$ be a prime and $K$ a finite extension of $\Q_p$ with residue field cardinality $q_K$. Write $W_K$ for the Weil group of $K$. For an algebraically closed field $\Omega$ of characteristic 0, denote by $\Rep_n(W_K)_\Omega$ (resp. $\Rep(\GL_n(K))_\Omega$) the set of isomorphism classes of $n$-dimensional Frobenius-semisimple Weil-Deligne representations of $W_K$ (resp. irreducible smooth representations of $\GL_n(K)$) on $k$-vector spaces.  For simplicity an element of $\Rep_n(W_K)$ will be called an ($n$-dimensional) WD-representation of $W_K$. Recall that such a representation is represented by $(V,\rho,N)$ where $V$ is an $n$-dimensional space over $\Omega$, $\rho:W_K\ra \GL_\Omega(V)$ is a representation such that $\rho(I_K)$ is finite and $\rho(w)$ is semisimple for every $w\in W_K$, and $N\in \End_\Omega(V)$ is a nilpotent operator such that $wNw^{-1}=|w|N$ where $|\cdot|:W_K\ra \R^\times_{>0}$ is the transport of the modulus character on $K^\times$ via class field theory. The local Langlands reciprocity map is a bijection
  $$\rec_K: \Rep(\GL_n(K))_\C\ra \Rep_n(W_K)_\C$$
  uniquely characterized by a list of properties, cf. \cite{HT01}. In particular $L(s,\pi)=L(s,\rec(\pi))$, $\varepsilon(s,\pi,\psi)=\varepsilon(\rec(\pi),\psi)$ for any nontrivial additive character $\psi:F\ra \C^\times$ (and a fixed Haar measure on $F$), and we also have an equality of conductors $f(\pi)=f(\rec_K(\pi))$. Here the local $L$ and $\varepsilon$ factors as well as conductors are independently defined on the left and right hand sides. Here we will only recall the definition of the conductor and $L$-factor for WD-representations, which is due to Grothendieck, leaving the rest of definitions and further references to \cite{Tat79} and \cite{Tay04}. For $(V,\rho,N)\in \Rep_{n}(W_{K})_{\Omega}$ the conductor is given by
    $$f(V):=\dim (V/V^{I_{K}}\cap \ker N) + \int_0^\infty \dim V/ V^{I^u_{K}} du,$$
  where $I^u_{K}$ is the upper numbering filtration on the inertia group $I_{K}$. Now let $\mathrm{Frob}_K$ denote the geometric Frobenius in $W_K/I_K$. The local $L$-factor is defined to be
  $$L(s,V):=\det (1-\Frob_K q_K^{-s}|V^{I_K}\cap \ker N)^{-1}$$
  so that we have the equality $L(s,\pi)=\det (1-\Frob_K q_K^{-s}|\rec(\pi)^{I_K}\cap \ker N)^{-1}$ for $\pi\in  \Rep(\GL_n(K))_\C$.

  Now fix a field isomorphism $\iota:\ol{\Q}_l\simeq \C$ and let $\rho:\Gal(\ol{F}/F)\ra \GL_n(\ol{\Q}_l)$ be a continuous semisimple Galois representation which is unramified at almost all primes and potentially semistable (equivalently de Rham) at places of $F$ above $l$. Such a $\rho$ is to be called \emph{algebraic}. At each finite place $v$ of $F$, there is a functor $\WD_v$ from continuous representations of $\Gal(\ol{F}_v/F_v)\ra \GL_n(\ol{\Q}_l)$ (assumed potentially semistable if $v|l$) to WD-representations of $W_{F_v}$. The construction of $\WD_v$ relies on Grothendieck's monodromy theorem when $v\nmid l$ and Fontaine's work in $l$-adic Hodge theory if $v|l$.

  The (global) conductor for $\rho$ is $\prod_v\mathfrak{p}_v^{f_v}$ where $\mathfrak{p}_v$ is the prime ideal of $\mathcal{O}_F$ corresponding to $v$, and $f_v=f(\rho|_{\Gal(\ol{F}_v/F_v)})$.
  To $\rho$ is associated a product function in a complex variable $s$, which is a priori formal infinite product:
  $$L(s,\rho):=\prod_{v:\textrm{finite}}L_v(s,\rho),\qquad L_v(s,\rho):= L(s,\iota \WD_v(\rho|_{\Gal(\ol{F}_v/F_v)})).$$
  When $\rho$ arises as a subquotient in the $l$-adic cohomology of an algebraic variety over $F$, one can apply Deligne's purity theorem to show that $L(s,\rho)$ converges absolutely for $\textrm{Re}(s)\gg 1$ (with often explicit lower bound). Further there is a recipe for the archimedean factor $L_\infty(s,\rho)$ in terms of Hodge-Tate weights of $\rho$ at places above $l$. (See the definition of $\Gamma(R,s)$ in \cite[\S2]{Tay04}, taking $R$ to be the induced representation of $\rho$ from $\Gal(\ol{F}/F)$ to $\Gal(\ol{F}/\Q)$.) This leads to a completed $L$-function
  $$\Lambda(s,\rho):=L(s,\rho)L_\infty(s,\rho).$$

  In the main body of the paper we were interested in the $L$-functions for Galois representations arising from varieties.
  Let $X$ be a smooth projective variety over $\Q$, so $X$ has good reduction modulo $p$ for all but finitely many primes $p$. Then a reciprocity law for $X$ on a concrete level would be a description of the number of points of $X$ in $\F_p$ (and its finite extensions) in terms of automorphic data at $p$ (i.e., local invariants at $p$ of several automorphic representations of general linear groups) as $p$ runs over the set of primes with good reduction, cf. \cite{Lan76}. This may be thought of as a non-abelian reciprocity law generalizing the Artin reciprocity law in class field theory as well as an observation about elliptic modular curves by Eichler-Shimura.
  Now we say that $\rho:\Gal(\ol{\Q}/\Q)\ra \GL_n(\ol{\Q}_l)$ \emph{comes from geometry} if
\begin{itemize}
  \item $\rho$ is unramified away from finitely many primes,
  \item there exists a finite collection of smooth projective varieties $X_i$ and integers $d_i,m_i\in \Z$ (indexed by $i\in I$) such that $\rho$ appears as a subquotient of  $$\bigoplus_{i\in I} H^{d_i}_{\et}(X\times_\Q \ol{\Q},\ol{\Q}_l)(m_i).$$
\end{itemize}
  As usual  $(m_i)$ denotes the Tate twist. One can speak of the obvious analogue with $\Q$ replaced by any finite extension $F$ over $\Q$.
  In the language of $L$-functions the following conjecture presents a precise form of the reciprocity law as above.

\begin{conjecture}\label{conj:modularity}
  Let $\iota:\ol{\Q}_l\simeq \C$ be an isomorphism. If $\rho:\Gal(\ol{F}/F)\ra \GL_n(\ol{\Q}_l)$ comes from geometry then $L(s,\rho)$ is automorphic, namely there exists an isobaric automorphic representation $\Pi$ of $\GL_n(\A_F)$ such that $L_v(s,\Pi)=L_v(s,\rho)$ at every finite place $v$ and $v=\infty$ (so that $L(s,\Pi)=L(s,\rho)$ and $\Lambda(s,\Pi)=\Lambda(s,\rho)$).
\end{conjecture}

  The Hasse-Weil conjecture predicts that $L(s,\rho)$ should have nice analytic properties such as analytic continuation, functional equation, and boundedness in vertical strips. If we believe in the Hasse-Weil conjecture, the converse theorem (discovered by Weil and then developed notably by Piatetskii-Shapiro and Cogdell) gives us a good reason to also believe that Conjecture \ref{conj:modularity} is true.

  The conjecture begs two natural questions, namely a useful characterization of $\rho$ coming geometry and a description of $\Pi$ that arise from such $\rho$. The conjectural answers have been provided by Fontaine-Mazur and Clozel, respectively. Indeed a conjecture by Fontaine-Mazur asserts that a continuous semisimple $l$-adic representation $\rho$ comes from geometry if and only if it is algebraic. Following Clozel a cuspidal automorphic representation $\Pi$ of $\GL_n(\A_F)$ is said to be \emph{L-algebraic} if, roughly speaking, the $L$-parameters for $\Pi$ at infinite places consist of algebraic characters in a suitable sense (see \cite{BG} for the definition; this differs from \cite{Clo90} in that no adjustment by the $\frac{n-1}{2}$-th power is made, cf. comments below Conjecture \ref{conj:correspondence}). An isobaric sum of cuspidal representations $\boxplus_{i=1}^r\Pi_i$ is algebraic if every $\Pi_i$ is algebraic.
   Then we can reformulate Conjecture \ref{conj:modularity} as one about the existence of the global Langlands correspondence preserving $L$-functions:

\begin{conjecture}\label{conj:correspondence}
  Fix $\iota$ as above. Then there exists a bijection $\Pi\leftrightarrow\rho$ between the set of L-algebraic isobaric automorphic representations of $\GL_n(\A_F)$ and the set of algebraic $n$-dimensional semisimple $l$-adic representations of $\Gal(\ol F/F)$ (up to isomorphism) such that the local $L$-factors are the same, so that $L(s,\Pi)=L(s,\rho)$ and $\Lambda(s,\Pi)=\Lambda(s,\rho)$.
\end{conjecture}

\begin{remark}\label{remark:correspondence} The strong multiplicity one theorem and the Chebotarev density theorem imply that if there is a correspondence $\Pi\leftrightarrow\rho$ as above then it should be a bijective correspondence and unique (but it does depend on the choice of $\iota$). It is expected that the set of cuspidal $\Pi$ maps onto the set of irreducible $\rho$.  A stronger property, often referred to as the local-global compatibility, is believed to be true at finite places $v$: it says that $\rec_{F_v}(\Pi_v)=\iota \WD(\rho|_{\Gal(\ol{F}_v/F_v)})$. (This is stronger only at ramified places.) In particular it should be true that $\rho$ and $\Pi$ have the same conductor (at finite places). Since we are concerned with unitary duals, we have adopted the unitary normalization for the Langlands correspondence and algebraicity. For arithmetic considerations it is customary to twist $\Pi$ by the $\frac{1-n}{2}$-th power of the modulus character in the conjecture. If so, one should replace ``L-algebraic'' by ``C-algebraic'', cf. \cite{BG}. 
\end{remark}

It is worth noting that Conjecture \ref{conj:modularity} suffices for our purpose in discussing geometric families. An important part of the Langlands program has been to confirm Conjecture \ref{conj:modularity} when $\rho$ is the $l$-adic cohomology of a Shimura variety (in any degree), which in turn led to many instances of the map $\Pi\mapsto \rho$ in Conjecture \ref{conj:correspondence}. Another remarkable result toward the conjectures is the modularity of elliptic curves over $\Q$ due to Wiles and Breuil-Conrad-Diamond-Taylor, who identified $L(s,\rho)$ with the $L$-function of a weight 2 modular form when $\rho$ is the \'etale $H^1$ of an elliptic curve over $\Q$. Recent developments include modularity lifting and potential modularity theorems. As we have no capacity to make a long list of all known cases of either Conjecture \ref{conj:modularity} or \ref{conj:correspondence}, we mention survey articles \cite{Tay04} and \cite{Har10} for the reader to begin reading about progress until 2009.

We close the discussion with a comment on the unitarity of local components and the issue of correct twist, cf. Remark (iv) below the definition of geometric families in Section \ref{section1}. Consider the automorphic representation $\Pi$ corresponding via the above conjectures to $\rho=H^d_{\rm et}(X\times_\Q \ol{\Q},\ol{\Q}_l)$ for a smooth proper variety $X$ over $\Q$ (which is not necessarily geometrically connected). Set $\Pi':=\Pi\otimes |\det|^{d/2}$. If $X$ has good reduction modulo a prime $p$ then the geometric Frobenius acts on the $H^d$-cohomology with absolute values $p^{d/2}$ under any choice of $\iota$. (This is Deligne's theorem on the Weil Conjectures if $p\neq l$. The argument extends to $p=l$ by work of Katz-Messing.) Hence the twist the Satake parameters of $\Pi'_p$ have absolute value 1, so $\Pi'_p$ is unitary. In general when $X$ has bad reduction modulo $p$, the unitarity of $\Pi'_p$ can be deduced from the weight-monodromy conjecture in mixed characteristic (as stated in \cite{Sai03}). Despite recent progress, cf. \cite{Sch12}, the latter conjecture is still open. What we said of $\rho$ should remain true when $\rho$ is a subquotient of $H^d_{\rm et}(X\times_\Q \ol{\Q},\ol{\Q}_l)$.

\section{Non-criticality of the central value for orthogonal representations}\label{sec:non-criticality}

    Deligne (\cite{Deligne:valeurs}) made a conjecture on special values of motivic $L$-functions. For a given $L$-function there is a set of the so-called \emph{critical} values of $s$ to which his conjecture applies. For our purpose we take on faith a motivic version of Conjecture \ref{conj:correspondence} (cf. \cite[\S6]{Lan12} and Remark \ref{remark:correspondence} above) on the existence of a bijection between absolutely irreducible pure motives $M$ of rank $n$ over $\Q$ and cuspidal C-algebraic automorphic representations $\pi$ of $\GL_n(\A)$ such that $$L(s+\frac{n-1}{2},M)=L(s,\pi).$$
Thereby Deligne's conjecture translates to a conjecture on automorphic $L$-functions. We copy the definition of $s$ being critical from the motivic side to the automorphic side in the obvious way. We are particularly interested in the question of whether the central value $s=1/2$ is critical for a cuspidal automorphic $L$-function which is unitarily normalized (for this a twist by a suitable power of the modulus character may be needed). The goal of appendix is to show

\begin{proposition} Suppose that a cuspidal automorphic representation $\pi$ of $\GL_n(\A)$ is
 \begin{enumerate}
   \item orthogonal (i.e., $\pi$ is self-dual and $L(s,\pi,\Sym^2)$ has a pole) and
   \item regular and C-algebraic. 
 \end{enumerate}
  Then $s=1/2$ is not critical for $L(s,\pi)$.
\end{proposition}

 The statement, in particular the definition of criticality, is unconditional in that no unproven assertions need to be assumed. However the proof is conditional on Conjecture \ref{conj:correspondence} as well as various conjectures around motives that are supposed to be true (see section 1 of \cite{Deligne:valeurs} for the latter). We freely assume them below.

\begin{proof}

 There should be a pure irreducible rank $n$ motive $M$ over $\Q$ corresponding to $\pi$. We follow the conventional normalization so that the weight of $M$ is $w=n-1$. (Note that the second assumption on $\pi$ implies that $M$ has Hodge numbers 0 or 1. In the Hodge realization the dimension of $M^{p,q}$ is at most one, and zero if $p+q\neq n-1$.) Since $\pi$ is self-dual, $M$ is self-dual up to twist. More precisely there is a perfect pairing $$  M \otimes M \rightarrow \Q(1-n)$$
  where $\Q(1-n)$ is the $(1-n)$-th power of the Tate motive.

  The center of symmetry for $L(s,M)$, the $L$-function associated to $M$, is at $s = (1+w)/2=n/2$. The necessary condition (which may not be sufficient) for it to be critical is that $n/2\in \Z$, namely that $n$ is even (so $w$ is odd). Hence we may and will assume that $n$ is even.
   Now consider the $l$-adic realization
  $$M_l \otimes M_l \rightarrow \Q_l(1-n),$$
where $M_l$ is now an irreducible $l$-adic representation of $\Gal(\ol{\Q}/\Q)$. By a result of Bellaiche-Chenevier's (\cite{BC11}) the
sign of $M_l$ is equal to $(-1)^{n-1}=-1$, meaning that the above pairing on $M_l$ is symplectic. (To apply their result we need both assumptions (1) and (2) on $\pi$.) Translating back to the automorphic side we deduce that $\pi$ is also symplectic. We have shown that if $s=1/2$ is critical then $\pi$ is symplectic, completing the proof.

\end{proof}


\begin{example}
  When $n=1$ and $\pi$ corresponds to a Dirichlet character $\chi$, it is well known that the central value $s=1/2$ for $L(s,\chi)$ is not critical. In this case $\pi$ is clearly orthogonal and the proposition applies. 
\end{example}

\begin{example}
  Consider the case of $n=2$ where $\pi$ corresponds to weight $k$ cuspforms ($k\ge 1$). Since we are concerned with self-dual representations, we normalize the correspondence such that $\pi$ is self-dual. Then $\pi$ is regular algebraic if and only if $k$ is even. (To deal with odd weight forms, one could twist $\pi$ by a half-power of the modulus character, but then $\pi$ would be self-dual only up to a twist.) In case $k$ is even, we associate to $\pi$ a pure motive $M$ of rank 2 and weight 1 such that $\dim
M^{1-k/2,k/2} = \dim M^{1-k/2,k/2} = 1$. 
It is equipped with a symplectic pairing $M \times M \rightarrow \Q(-1)$.


\end{example}

\subsection*{Acknowledgments} Thanks to Manjul Bhargava, Daniel Bump, Brian Conrey, Dorian Goldfeld,  Henryk Iwaniec, Philippe Michel, Emmanuel Kowalski, Erez Lapid, Gopal Prasad, Zeev Rudnick, Jean-Pierre Serre, Arul Shankar, Anders S\"odergren, Kannan Soundararajan, Akshay Venkatesh, Jun Yu, and especially Nicholas Katz, for discussions and insights on various aspects of the paper. S.W.S. is grateful to Princeton University and the Institute for Advanced Study for their hospitality during several short visits. He acknowledges partial supports from NSF grant DMS-1162250 and a Sloan Fellowship. N.T. acknowledges partial support from NSF grant DMS-1200684.

\def\cprime{$'$} \def\cprime{$'$} \def\cprime{$'$} \def\cprime{$'$}
  \def\cprime{$'$} \def\cprime{$'$} \def\cprime{$'$} \def\cprime{$'$}
  \def\cprime{$'$} \def\cprime{$'$}

\end{document}